\documentstyle[11pt]{amsart}

\date{}

\newcommand{\commentout}[1]{}
\newcommand{\comment}[1]{{\bf [***COMMENT:#1***]}}
\newcommand{\new}{\mbox{\tt\char'134begin\{New\}}}   
\def\endnew{\mbox{\tt\char'134end\{New\}}}  
\newcommand{\mnote}[1]{\marginpar{\footnotesize\em #1}}
\newcommand{\note}[1]{}
\newcommand{\DELETENOTES}{
\newcommand{\lars}{}
\newcommand{\ralph}{*RH*\ }
\renewcommand{\commentout}[1]{}
\renewcommand{\comment}[1]{}
\renewcommand{\new}{}
\renewcommand{\endnew}{}
\renewcommand{\mnote}[1]{}
\renewcommand{\note}[1]{}
}
\DELETENOTES
\newcommand{\ac}{\overline{\nabla}}	
\newcommand{\con}{\nabla}	   	
\newcommand{\am}{\overline{\phantom{N}}\!\!\!\!\!\! M}
\newcommand{\la}{\langle}	   	
\newcommand{\ra}{\rangle}          	
\newcommand{\oc}{^{\bot}}          	
\newcommand{\e}{\varepsilon}       	
\newcommand{\nc}{\nabla^{\bot}}	   	
\newcommand{\ar}{\,\overline{\phantom{I}}\!\!\!\! R}
\newcommand{\aS}{\,\overline{\phantom{I}}\!\!\!\! S}
\newcommand{\ch}{{\cal H}}         	
\newcommand{\cv}{{\cal V}}	   	
\newcommand{\r}{{\Bbb R}}		
\newcommand{\trace}{\mathop{\mathrm{trace}}}	
\newcommand{\Ric}{\mathop{\mathrm{Ric}}}	
\newcommand{\Scal}{\mathop{\mathrm{Scal}}}	
\newcommand{\ame}{\bar{g}}	   	
\newcommand{\bargw}{\bar g_w}           
\newcommand{\T}{\partial_t}		
\newcommand{\w}{\wedge}			
\newcommand{\W}{\bigwedge}		
\newcommand{\nor}{\overline{\eta}}	

\newcommand{\h}{\bold{H}}
\newcommand{\hm}{g_{\bold{H}}}
\newcommand{\rr}{\overline{r}}		
\newcommand{\elm}{g_{\bold{S}}}
\newcommand{\els}{\bold{S}}
\newcommand{\whp}{\widehat{\phi}}	
\newcommand{\cl}{\mathop{\mathrm{Cl}}}	
\renewcommand{\phi}{\varphi}		
\newcommand{\eps}{\varepsilon}		
\newcommand{\half}{\frac{1}{2}}
\newcommand{\cd}{,\ldots,}              

\newcommand{\aRic}{\overline{\mathop{\mathrm{Ric}}}}

\renewcommand{\(}{\left(}
\renewcommand{\)}{\right)}
\renewcommand{\[}{\left[}
\renewcommand{\]}{\right]}

\newcommand{\done}{\qed}

\renewcommand{\mod}{\r^n_k(K_0)}
\newcommand{\mme}{g_{\mod}}
\newcommand{\stack}[2]{\begin{array}{c} #1 \\ #2 \end{array} } 
\newcommand{\vs}{\vphantom{\dfrac1{\int}}}

\theoremstyle{plain}
\newtheorem{thm}{Theorem}[section]
\newtheorem{lemma}[thm]{Lemma}
\newtheorem{prop}[thm]{Proposition}
\newtheorem{cor}[thm]{Corollary}
\newtheorem{remark}[thm]{Remark}
\newtheorem{definition}[thm]{Definition}

\numberwithin{equation}{section}   

\title[Comparison and Rigidity Theorems]{Comparison and Rigidity Theorems in Semi-Riemannian Geometry}
\author[Lars Andersson]{Lars Andersson\footnotemark{$^*$}}
\thanks{\footnotemark{$^*$} Supported in part by NFR, contract no. F-FU 4873-307.} 
\address{Department of Mathematics \\
Royal Institute of Technology \\
S-100 44 Stockholm, Sweden}
\email{larsa\char'100math.kth.se}
\author[Ralph Howard]{Ralph Howard\footnotemark{$^\dagger$}}
\thanks{\footnotemark{$^\dagger$} Supported in part by NFR, contract no. R-RA 4873-306 and the Swedish 
Academy of Sciences.}
\address{Department of Mathematics\\
University of South Carolina \\
Columbia, S.C. 29208, USA}
\email{howard\char'100math.sc.edu}
\date{June 4, 1996}
\begin{document}

\begin{abstract}
The comparison theory for the Riccati equation satisfied by the shape
operator of parallel hypersurfaces is generalized to semi--Riemannian
manifolds of arbitrary index, using one--sided bounds on the Riemann
tensor which in the Riemannian case correspond to one--sided bounds on
the sectional curvatures. Starting from 2--dimensional rigidity
results and using an inductive technique, a new class of gap--type
rigidity theorems is proved for semi--Riemannian manifolds of
arbitrary index, generalizing those first given by Gromov and
Greene--Wu. As applications we prove rigidity results for
semi--Riemannian manifolds with simply connected ends of constant
curvature.  
\end{abstract}

\maketitle

{\small 
\tableofcontents
}

\section{Introduction}\label{sec:intro}

In Riemannian geometry, the comparison results in terms of sectional
curvature of Rauch, Toponogov, Morse--Schoenberg and others
(cf.~\cite{Cheeger-Ebin}) are basic tools used to prove results such
as the sphere theorem, the Bonnet--Myers theorem, and the maximal
diameter theorem of Toponogov.
More recently, comparison theorems in terms of the Ricci curvature
such as the Bishop--Gromov volume comparison theorem have played an
important role leading to such results as the Chen maximal diameter
theorem, see the wonderful survey article by
Karcher~\cite{Karcher:comp}.

In Lorentzian geometry and semi-Riemannian geometry, on the other
hand, it is well known (cf.~\cite{Harris:triangle} and
Section~\ref{sec:examples} for a discussion) that the assumption of
even a one--sided bound on sectional curvature defined as
$$
\frac{\la R(X,Y)Y,X\ra}{\la X,X\ra\la Y,Y\ra-\la
X,Y\ra^2},
$$
implies that the space has constant sectional curvature.
Therefore such bounds are not interesting.

In the general semi-Riemannian setting the natural replacement of a
one side bound on the sectional curvature is
\begin{definition}\sl
\label{def:curvbound}
We will say that $R \geq K_0$ or $R \leq K_0$ if and only if for all 
$X,Y$, the inequalities
\begin{eqnarray}
        &\la R(X,Y)Y,X\ra 
	\ge K_0(\la X,X\ra \la Y,Y\ra -\la X,Y\ra^2) &\label{RgeK0} ,  \\
&\mbox{ or } & \nonumber\\
        &\la R(X,Y)Y,X\ra \le K_0(\la X,X\ra \la Y,Y\ra -\la X,Y\ra^2) & 
\label{RleK0}
\end{eqnarray}
hold.~\done 
\end{definition}

In case $g$ is positive definite these conditions are equivalent to $g$
having sectional curvature bounded from below or above by $K_0$.
Assuming bounds of the form given in this definition we are able to
prove ``gap type" rigidity results of the type first proved by Greene
and Wu \cite{Greene-Wu:gap-theorems} and Gromov
\cite{ballmann:gromov:schroeder}.  Our main example of such a result
is

\begin{thm}\sl \label{introthm} 
Let $(M,g)$ be a geodesically complete semi-Riemannian manifold of
dimension $n\ge3$ and index~$k$ with
curvature satisfying one of the two inequalities $R \geq 0$
or $R\le 0$  and assume  $(M,g)$ has an end $E$ with
$R\equiv 0$ on $E$ and $\pi_1(E)$ finite.  Then $(M,g)$ is isometric
to the flat model space $(\r^n_k,g_0)$.
\end{thm} 

In the Lorentzian and Riemannian case we also have rigidity results
for comparison with constant, but non-zero, curvature models, cf.~Section~\ref{sec:rigid}.  Even in the Riemannian case this leads to a
new rigidity result.
\begin{thm}\sl\label{thm:IntroRiemThm}
Let $(M,g)$ be a complete Riemannian manifold of dimension
$n\ge3$ with sectional
curvatures $\le1$.  Let $B\subset S^n\setminus S^{n-1}$ be a closed set with
$S^n\setminus B$ connected and let $\phi :S^n\setminus B\to M$
be a local isometry.  Then $\phi$ extends to a surjective local
isometry $\widehat{\phi}:S^n\to M$. Therefore $(M,g)$ is a quotient of
$S^n$ by a finite fixed point free group of isometries.
\end{thm}
This result differs from most rigidity theorems for the sphere in
that an upper bound, rather than a lower bound, is given on the
curvature.

The proofs of the gap theorems rest on two main ideas. The first is an
extension of the comparison theory for matrix Riccati equations to the
semi-definite setting, cf.~Section~\ref{sec:basic}. The second is an
inductive argument involving the use of foliations by parallel
hypersurfaces, with base case given by some 2-dimensional rigidity
results, see Section~\ref{subsec:two-dim}. The two are related by the
fact that the second fundamental form of a parallel foliation
satisfies the Riccati equation.

\subsection{Background}
One important class of results in Riemannian geometry are the ``gap type'' 
rigidity theorems. Let us mention a few examples of model gap
results. We do not try to state these results under the weakest
hypotheses. Let $(M,g)$ be a complete
connected Riemannian manifold, isometric to Euclidean $\r^n$ 
in the complement of a compact set with $n\ge3$.
\begin{itemize}
\item {\bf Scalar curvature:} 
Assume that $(M,g)$ is spin and has non--negative scalar curvature.
Then $(M,g)$ is isometric to Euclidean $\r^n$. This follows from the
Witten argument for the positive mass theorem. 
\item {\bf Ricci curvature} Assume that $(M,g)$ has non--negative 
Ricci curvature.
Then the conclusion follows from the Bishop--Gromov volume 
comparison theorem. 
\item {\bf Sectional curvature} If $(M,g)$ has either 
non--negative or non--positive
sectional curvature, then theorems 
of Greene and Wu~\cite{Greene-Wu:gap-theorems} (for
non-negative sectional curvature) and  Kasue and 
Sugahara~\cite{Kasue-Sugahara:Gap} (for non-positive 
sectional curvature)
imply $(M,g)$ is isometric to $\r^n$.  
\end{itemize}
In each of the above cases there are versions where the assumption that
$(M,g)$ is isometric to the standard space in the complement of a compact
set is replaced by asymptotic
assumptions. There are also versions where the Euclidean space is replaced
by constant curvature or homogeneous spaces. Cf.~Sections
\ref{sec:rigid-flat}, \ref{sec:rigid}, \ref{sub:ends}, and \ref{sub:quot}  
for further discussion. 

In the Lorentzian case, the most important rigidity type result is of
course the positive mass theorem of Schoen--Yau and Witten
(cf.~\cite{schoen:yau:bondi,witten:mass,parker:taubes}).  Witten's
proof has been generalized to the case of the Bondi mass by various
authors~\cite{horowitz:perry}, \cite{reula:tod},
\cite{ludvigsen:vickers}.

In the Riemannian case, the dependence of the scalar curvature on $g$ can be
viewed as a scalar elliptic operator acting on the conformal factor, and the
Ricci tensor can be viewed as a quasilinear elliptic system in terms of the
metric tensor. This makes the form of the gap type rigidity theorems natural,
in view of the maximum principle and unique continuation results which hold
for elliptic equations.

In the Lorentzian case, the Ricci and Einstein tensors form hyperbolic
systems in terms of the metric tensor and the positive mass theorems is
closely related to the elliptic constraints induced via the Gauss and
Gauss--Codazzi equations on the extrinsic geometry of spacelike
hypersurfaces. Further, in this case, causality is a powerful organizing
principle in the geometry which shows up for example in the conservation
theorem~\cite[\S 4.3]{Hawking-Ellis} and Theorem~\ref{thm:EinRigid} below.

 From the PDE point of view, putting conditions on the whole Riemann tensor
constitutes an over determined system, regardless of the signature.
More precisely in normal coordinates centered at a point of a
semi-Riemannian manifold it is possible to express all the
derivatives of order at least two of the metric in terms of the
curvature tensor and its derivatives.  Thus the operator that takes a
metric to its curvature tensor behaves very much like an elliptic
operator in that it is possible to bound the second and higher
derivatives of the unknown metric in terms of the curvature tensor.
Therefore, one expects bounds on the Riemann tensor to give gap type
rigidity theorems regardless of the signature, as is shown to be the
case in the present paper.

In the Lorentzian case, the one-sided curvature bounds 
used in this paper implies conditions on the timelike sectional curvatures
in the sense of \cite{Beem-Erlich:Lorentz2}.
For example, $R \geq 0$, 
implies nonpositive timelike sectional curvatures in the sense of
\cite{Beem-Erlich:Lorentz2}.
The  condition $R \geq 0$ also implies the strong energy condition, i.e.
that $\Ric(X,X) \geq 0$ for timelike vectors $X$. In this situation, the
Lorentzian splitting theorems of Galloway, Eschenburg, Newman and 
others~\cite{galloway:review} 
would apply if we add the condition of existence of a timelike ray in the
future of a maximal hypersurface, cf.\ 
\cite[Theorem C]{galloway:1989}.
This would then easily imply the Lorentzian case of Theorem~\ref{flatrigid}. 
Therefore it is important that our
assumptions do not include the existence of a timelike ray.

Further, in view of the fact that $R \geq 0$ implies the strong energy
condition, under the additional assumption of global hyperbolicity, Theorem
\ref{flatrigid} would in the Lorentzian case follow easily from the 
singularity theorems in Lorentzian geometry, cf.~\cite{Hawking-Ellis}. 
Therefore we also stress that our conditions do not include global
hyperbolicity.

Let us emphasize that the conditions of timelike geodesic 
completeness and global hyperbolicity are not at all closely related. 
As an illustration of this fact, let us mention that anti--de~Sitter space
is timelike geodesically complete but {\em not} globally hyperbolic, while 
the ``strip" $\els_1^n(-1)$,
cf.~Definition \ref{def:half-space}, which is a subset of anti--de~Sitter
space, is globally hyperbolic, but not 
timelike geodesically complete. 
The assumption of timelike geodesic completeness is thus natural in the
context of anti--de~Sitter space.

\subsection{Overview of this paper.}
Section~\ref{prelim} contains  the basic results and
formulas we will be using.  Almost all of this is well known, but is
included to fix notation and as the details of writing the Riemannian
and Lorentzian space forms as warped products is quite important to
many of our rigidity results.

Section~\ref{sec:basic} contains the generalization of the basic
comparison theory to the general semi-Riemannian setting.  Most of
the results are stated as results about systems of ordinary
differential equations.  

Section \ref{sec:rigid-flat} gives the rigidity results for flat
semi-Riemannian manifolds in terms of one sided bounds on the
curvature tensor.

Section~\ref{sec:space-forms} has rigidity results for the one sided
bounds on the sectional curvature when the model space has non-zero
constant sectional curvature.  The results in this section only apply
to Riemannian and Lorentzian manifolds.  

Section~\ref{sec:applications} contains various applications of the
rigidity results to manifolds that have an end of constant curvature.
This involves a classification of the ends of constant curvature and
finite fundamental groups and which is new in the case of indefinite
metrics.  There are also rigidity results for semi-Riemannian metrics
on quotients of space forms.   

\bigskip
\noindent
{\bf Acknowledgements:\ }  E. Calabi~\cite{Calabi:Sphere-Rigidity}
supplied us with Theorem~\ref{thm:calabi} which replaces a weaker 
ad hoc version.  Much of the work on this paper was done
while the second author (RH) was on sabbatical leave from the University
of South Carolina at the Royal Institute of Technology in Stockholm.
He would like to express his thanks to both schools.

\section{Preliminaries}\label{prelim}
\subsection{Basic Formulas and Results}

The {\em index\/} of a symmetric bilinear form $g(\,,)$ is the
maximal dimension of a subspace on which it is negative definite, or
equivalently, the number of negative eigenvalues.
If $\am$ is a smooth manifold then a {\em
semi-Riemannian metric\/} $\ame$ on $\am$ is a smooth $(0,2)$ tensor
which is everywhere non-degenerate.  The index of $\ame$ is locally
constant so if $\am$ is connected then the index~$k$ of $\ame$ is
well defined.  
Thus the semi-Riemannian manifolds of index~0 are the
Riemannian manifolds and the ones of index~$1$ are the {\em
Lorentzian\/} manifolds.   Unlike the Riemannian case not every
manifold need have a semi-Riemannian metric of index~$k$.  A
necessary and sufficient condition for a smooth manifold $\am$ to have a
semi-Riemannian metric of index~$k$ is that the tangent bundle
$T(\am)$ have a rank~$k$ subbundle $E$.  For $\am$ compact
this is proven in \cite[Theorem~40.11]{Steenrod:bundles}, but the
extension to noncompact manifolds is easy. In particular, every noncompact
manifold has a Lorentzian metric.

Let $\am$  be an $n$~dimensional semi-Riemannian manifold 
with semi-Riemannian metric $\ame(\,,)=\la\,,\ra$ 
of index~$k$.  An immersed
submanifold $M$ of $\am$ of dimension~$m$ is called {\it non-degenerate} if 
the restriction of the semi-Riemannian metric $\la\,,\ra$ to $T(M)$ is 
non-degenerate at each point of $M$.  The normal bundle to $M$ in $\am$ 
will be denoted by $T\oc(M)$.  This is the bundle over $M$ so that at 
each $x\in M$ the fiber $\{Y\in T(\am)_x:\la Y,X\ra=0\  
\mbox{for all}\ X\in T(M)_x\}$.
There is a direct sum decomposition 
$$
T(\am)_x=T(M)_x\oplus T\oc(M)_x   
$$
at a
point $x\in M$ if and only if the restriction of $\la\,,\ra$ to $T(M)_x$ is
non-degenerate.  This also equivalent to  the restriction of 
$\la\,,\ra$ to the normal bundle $T\oc(M)$ being non-degenerate.  
For the rest of this section we will assume that the induced metric on 
$M$ is non-degenerate.

Let $\ac$ be the metric connection on $\am^n$ and let $\con$ be the metric
connection on $M$.  These are related by
$$
\ac_XY=\con_XY+h(X,Y) ,
$$
where $X$, $Y$ are smooth vector fields tangent to $M$, and $h(\,,)$ is the
vector valued second fundamental form of $M$ in $\am$.  Thus $h(\,,)$ is a
symmetric bilinear map from $T(M)\times T(M)\to T\oc(M)$.  If $\xi$ is a
smooth section of $T\oc(M)$ over $M$ and $X$ is a smooth vector field
tangent to $M$ then 
$$
\ac_X\xi=\nc_X\xi-A_{\xi}X,
$$
where $\nc_X\xi\in T\oc(M)$ and $A_{\xi}X\in T(M)$.  Then $\nc$ is the
{\it connection in the normal bundle}, and $A_{\xi}$ is the {\it Weingarten
map} or {\it shape operator} in the normal direction $\xi$.  This is related
to the second fundamental form by 
$$
\la A_{\xi}X,Y\ra=\la h(X,Y),\xi\ra.
$$
Thus the symmetry of $h(\,,)$ implies that for any $\xi\in T\oc(M)$ the
linear map $A_{\xi}:T(M)\to T(M)$ is self-adjoint with respect to the inner
product $\la\,,\ra$.

Let $\ar$ be the curvature tensor of $\ac$ where our choice of sign on the
curvature tensor is 
$$
\ar(X,Y)Z=\ac_X\ac_YZ-\ac_Y\ac_XZ-\ac_{[X,Y]}Z .
$$
This is related to the curvature tensor of $\con$ by the Gauss
curvature equation
\begin{equation}
\la R(X,Y)Z,W\ra=\la \ar(X,Y)Z,W\ra
+\la h(X,W),h(Y,Z)\ra-\la h(X,Z),h(Y,W)\ra.
\end{equation}
When $M$ is a hypersurface with unit normal field $\xi$ this can be
rewritten in terms of the Weingarten map as 
\begin{equation}\label{gauss} 
\begin{aligned} 
\la R(X,Y)Z,W\ra=&\la \ar(X,Y)Z,W\ra\\
&\quad +\la \xi, \xi \ra \left( \la A_{\xi}X,W\ra\la A_{\xi}Y,Z\ra-
\la A_{\xi}X,Z\ra\la A_{\xi}Y,W\ra\right) . 
\end{aligned}
\end{equation}

The normal bundle has two unit sphere  bundles $S\oc_{+1}(M)$ and
$S\oc_{-1}(M)$. 
Letting  $\e$ be either  $\e=+1$ or $\e=-1$ they are
defined by 
$$
S\oc_{\e}(M)=\{u\in T\oc(M):\la u,u\ra=\e\}.
$$
Each of these is a bundle $\pi:S\oc_{\e}(M)\to M$
with fiber diffeomorphic to the ``sphere'' 
$\{u\in T\oc(M)_x:\la u,u\ra=\e\}$.  

Let $\pi:T\oc(M)\to M$ be the projection and for $v\in
T\oc(M)$ let $\cv_v:=\ker(\pi_{*v})$ be the vertical vectors at $v$.
There is a natural identification of $\cv_v$ with 
the subspace $T\oc(M)_{\pi x}$ of
$T(\am)_{\pi v}$.  In a vector space the tangent spaces $T(V)_x$ are
naturally identified with $V$. The fibers $\pi^{-1}[x]=T\oc(M)_x$ 
of $\pi$ are
vector spaces and $\cv_v=T\oc(M)_{\pi x}\subset T(\am)_{\pi v}$.  
Let $\ch_v$ be the space of horizontal vectors
at $v$. A curve $\gamma:(a,b)\to T\oc(M)$ is horizontal iff it
it is parallel along its projection $c=\pi\circ\gamma$.  That is it satisfies
$\con\oc_{c'(t)}\gamma=0$. A vector $X$ tangent to $T\oc(M)$ 
is horizontal iff it is tangent to a horizontal curve.   
The space $\ch_v$ can be identified
with $T(M)_{\pi V}$ by the map $\pi_*|_{\ch_v}$.  These identifications
combine to give $T(T\oc(M))_v=\ch_v\oplus\cv_v=T(M)_{\pi v}\oplus
T\oc(M)_{v}=T(\am)_{\pi v}$.  Under this identification for $u\in
S\oc_\e(M)$ the tangent space $T(S^\bot_\e)_u$ is identified with
$u\oc\subset T(M)_{\pi u}$.  

Let $\exp$ be the exponential map defined w.r.t. 
$\ac$ and define a function
$f_r:S\oc_\e\to \am$ by $M(r):=f_r[S_\e^\bot]$.  Then $M(r)$ is the
{\em tube of radius\/} $r$ (and sign $\e$) about $M$.  
Let $u\in T(\am)_x$ and let $\gamma(t)=\exp(tu)$ be
the geodesic that fits $u$.  Identify all the tangent spaces
$T(\am)_{\gamma(t)}$ along $\gamma$ with $T(\am)_x$ by parallel translation
along $\gamma$.  For each $t$ let $\ar_u(t):T(\am)_x\to T(\am)_x$ be the
linear map 
$$
\ar_u(t)X=\ar_{\gamma(t)}(X,u)u=\ar_{\gamma(t)}(X,\gamma'(t))\gamma'(t) ,
$$
where $\gamma'(t)$ and $u$ are identified as $\gamma'$ is parallel along
$\gamma$. 
For each $t$ the linear map $\ar_u(t)$ is self-adjoint ({\it i.e.\/} 
$\la\ar_u(t)X,Y\ra=\la X,\ar_u(t)Y\ra$).  Also for each $t$, $\ar_u(t)u=\ar(u,u)u=0$.  Thus
from elementary linear algebra the subspace $u\oc$ is invariant under
$\ar_u(t)$.

\begin{definition}\sl 
For $u\in S\oc_{\e}(M)$ let $P:u\oc\to T(M)_{\pi u}$ and 
$P\oc:u\oc\to T\oc(M)\cap u\oc$ be the orthogonal projections.  Define a
field of linear maps $F_u(t)$ along $\gamma(t)=\exp(tu)$ by the initial
value problem:
$$
F_u''(t)+\ar_u(t)F_u(t)=0,\quad F_u(0)=P,\quad F_u'(0)=-A_uP+P\oc,
$$
where ``prime'' denotes the covariant derivative by $\gamma'(t)$
({\em i.e.\/} $'=\ac_{\gamma'}$) and $A_u$ is the Weingarten map of 
$M$ in the direction $u$.~\done
\end{definition}

\begin{prop} With the identification of $T(S\oc_{\e}(M))$ with $u\oc$ given 
above, the derivative of $f_r$ at $u$ is
$$
{f_r}_{*u}=F_u(r).
$$
\end{prop}

\begin{pf} Just as in the Riemannian case.  Cf.~\cite{Heintze-Karcher}.
\end{pf}

If $F_u(r)$ is non-singular then the tube $M(r)$ is a smooth immersed
hypersurface at $f_r(u)=\exp(ru)$.  The outward normal at this point is 
\begin{equation}\label{etadef}
\eta(f_r(u))=\left.\frac{d}{dt}\exp(tu)\right|_{t=r}  .
\end{equation}
At these smooth points of $M(r)$  the  Weingarten map or 
of $M(r)$ with respect to this normal is the linear map 
$S_u(r):T(M(r))_{f_r(u)}\to T(M(r))_{f_r(u)}$ given by
$$
S_{u}(r)X=-\ac_X \eta.
$$
\begin{prop}\sl\label{2ndForm}
Identify $T(M(r))_{f_r(u)}$ with the subspace $u\oc$ of $T(M)_{\pi u}$ by
parallel translation along $t\mapsto \exp(tu)$.  Then the Weingarten map
is given by
\begin{equation}
S_u(r)=-F_u'(r)F_u(r)^{-1} .  \label{shape}
\end{equation}
At points where $F_u(t)$ is non-singular this satisfies the Matrix Riccati 
equation
\begin{equation}
S_u'(t)=S_u(t)^2+\ar_u(t)  \label{Ricciti}
\end{equation}
\end{prop}
\begin{pf} Follows from the last proposition by a calculation.  
Cf.~\cite{Heintze-Karcher}.
\end{pf}

The mean curvature of $M$ in the direction $\xi$ is 
\begin{equation}
\label{meancurv}
H_{\xi} = \frac{1}{n-1}\trace(A_{\xi})  .
\end{equation}
and we denote by $H_u(r)$ the mean curvature of $M(r)$ at $f_r(u)$.
If $A$ is an $(n-1)\times(n-1)$ symmetric matrix, then by the 
Cauchy-Schwartz inequality $(\trace(A))^2\le (n-1)\trace(A^2)$ with equality
if and only if $A=cI$ for some scalar $c$.  This leads to:
\begin{cor}\sl 
Assume that $\am$ is Riemannian or Lorentzian, and that $M$ is a space
like hypersurface in $\am$ (so that the induced metric on
$M(r)$ is Riemannian). Let $x\in M$, $\eta(x)$ the normal to $M$ at
$x$, $S(t)$ the Weingarten map of $M(t)$ along
$\gamma(t)=\exp(t\eta(x))$ and $H(t)=\frac1{n-1}$ the mean curvature.
\begin{equation}\label{tracericcati}
H'(t) \geq (n-1) H^2(t) + \trace(\ar_u(t))
\end{equation}
If equality holds at $t_0>0$ then $A_{\eta(x)}=S(0)=cI$ and
$\ar(t)=\overline{r}(t)I$ for $0\le t\le t_0$.~\done
\end{cor}
A fact that will be used repeatedly in the sequel is that along
$\gamma(t)=\exp(tu)$ 
\begin{equation}
\aRic (\gamma'(t),\gamma'(t))=\trace(\ar_u(t))	 .    \label{eq:ric-trace}
\end{equation}

We close this section with a result about when local isometries
between semi-Riemannian manifolds are surjective. 

\begin{prop}\sl \label{onto} Let $M$ be a geodesically complete
semi-Riemannian manifold and $f:M\to N$ a local isometry where $N$ is a
connected semi-Riemannian manifold of the same dimension as $M$.  Then  
$f$ is surjective and $N$ is also geodesically complete.
\end{prop}

\begin{pf} As $f$ is a local diffeomorphism the image $f[M]$ is
open.  Thus it is enough to show that $f[M]$ is closed.  Let $y$ be
a point in the closure of $f[M]$.  Then as the exponential map $\exp:TN_y\to
N$ maps a neighborhood of $0\in TN_y$ onto a neighborhood of $y$ in $N$
there is a geodesic segment $c:[0,\delta]\to N$ so that $c(0)=y$ and 
$c(\delta)\in f[M]$.  Using that $f$ is a local isometry we see there
is a geodesic $\gamma:\r\to M$ so that
$(f\circ\gamma)'(0)=c'(\delta)$.  But then $f\circ\gamma:\r\to N$
is a geodesic and by the uniqueness of solutions to the geodesic equation
we have $(f\circ\gamma)(-\delta)=c(0)=y$.  Thus $y=f(\gamma(-\delta))\in
f[M]$.  Therefore $f[M]$ contains all its limit points and thus it is
closed.  This completes the proof that $f$ is surjective.

Let $c:(-\delta,\delta)\to N$ be a small geodesic segment.  Then as $f$ is
surjective and $M$ geodesically complete, there is a geodesic
$\gamma:\r\to M$ so that $f\circ\gamma$ extends $c$ to all of 
$\r$.  Thus $N$ is geodesically complete.
\end{pf}

\subsection{Warped Products and Model Spaces}
\label{sec:warpmodel}

\newcommand{\ga}{\alpha}                
\newcommand{\gb}{\beta}                 
\newcommand{\g}{\gamma}
\newcommand{\tb}[2]{$\stackrel{\mbox{#1}}{\mbox{#2}}$}
\newcommand{\vvv}{\vphantom{$\dfrac11$}}
\newcommand{\vvvv}{\vphantom{$\dfrac{\int}1$}}
\newcommand{\vvvd}{\vphantom{$\dfrac1{\int}$}}

Let $K_0 \in \r$ and let $k$ be an integer, $0 \leq k \leq n$.
We will denote by $\r^n_k(K_0)$ the simply connected 
model space of sectional
curvature $K_0$ and index~$k$ and we will use $g_{\r^n_k(K_0)}$ to
denote its metric. 
In the Riemannian case ($k = 0$) the model
spaces are either Euclidean space, spheres or hyperbolic spaces. In
the Lorentzian case,
the non-flat model spaces are de~Sitter ($K_0>0$) 
and anti--de~Sitter spaces ($K_0 < 0$).

Let $n\geq 3$. Then $\r^n_1(1)$ is the 
hypersurface of $\r^{n+1}_{1}$ given by 
\begin{equation}\label{k-hyp}
- x_1^2 + \sum_{i=2}^{n+1} x_i^2 = 1 ,
\end{equation}
with its induced metric. This is known as de~Sitter space and is
diffeomorphic to $S^{n-1}\times \r$. Similarly, the hypersurface 
\begin{equation}\label{k-hyp2}
-x_1^2 -x_2^2 + \sum_{i=3}^{n+1} x_i^2 = 1
\end{equation} 
of $\r^n_2$ with its induced metric
has sectional curvature~$-1$ and index~$1$ and is diffeomorphic to
$S^1\times\r^{n-1}$. The universal 
cover of this space is known as anti--de~Sitter space and will be denoted by
$\r^n_1(-1)$.

Some of the model spaces we will consider are best
represented as warped products of a hypersurface with the real line.  
For completeness we include the basic formulas in the geometry of warped
products.

Let $M^{n-1}$ be a semi-Riemannian manifold with metric $g(\,,)$,
and $(a,b)$ be an open interval in $\r$, and $w$ a positive real valued
function on $(a,b)$. Define a semi-Riemannian
metric $\bargw(\,, )$ on $\am^n=M^{n-1}\times (a,b)$ by 
\begin{equation}
\bargw=w(t)^2g+\e (dt)^2  ,              \label{eq:warp-metric}
\end{equation}
where $\e=+1$ or $\e=-1$.  A vector field $X$ on $M^{n-1}$ can also be
viewed as a vector field on $\am$ in the natural way.  Let $\con$ be the
semi-Riemannian connection of $g$ and $\ac$ the semi-Riemannian connection
of $\bargw$.   If $X$, $Y$ are vector fields on $M$ then the connections 
$\con$ and $\ac$ are related by
$$
\ac_XY=\con_XY-\frac{\e w'\bargw(X,Y)}{w}\T,\quad
\ac_X\T=\ac_{\T}X=\frac{w'}{w}X,\quad \ac_{\T}\T=0 .
$$
The vector field $\eta:=\T$ is a unit normal vector to the hypersurfaces
$M^{n-1}\times\{c\}$.  From the above formulas  the second
fundamental form $h$ and the Weingarten map $S$ of these hypersurfaces are
given by 
\begin{equation}
h(X,Y)=-\frac{\e w'}{w}g_w(X,Y)\T,\quad SX=-\frac{w'}{w}X.  \label{secondf}
\end{equation}
The curvature $R$ of $M$ and the curvature $\ar$ of $\am$ are related by
$$
\ar(X,\T)\T=-\frac{w''}{w}X ,
$$
so that 
$$
\bargw(\ar(X,\T)\T,X)=-\frac{w''}{w}\bargw(X,X)
=-\frac{\e w''}{w}\bargw(X,X)\bargw(\T,\T) .
$$
Therefore 
$$
\Ric(\T,\T)=-(n-1)\frac{w''(t)}{w(t)} .
$$
Using the Gauss curvature equation (\ref{gauss}) there is also the relation
\begin{equation}\label{turd} 
\begin{aligned}
\bargw(R(X,Y)Y,X)=&\bargw(\ar(X,Y)Y,X) \\ & \quad  + \e\(\frac{w'}{w}\)^2
	(\bargw(X,X)\bargw(Y,Y)-\bargw(X,Y)^2) . 
\end{aligned}
\end{equation}

If $(M,g)$ has constant sectional curvature $K_0$ then
$$
g(R(X,Y)Y,X)=K_0(g(X,X)g(Y,Y)-g(X,Y)^2) , 
$$ 
so that $\bargw(R(X,Y)Y,X)=(K_0/w^2)(\bargw(X,X)\bargw(Y,Y)-\bargw(X,Y)^2)$.
Using this in  equation (\ref{turd}) leads to 
$$
\bargw(\ar(X,Y)Y,X)=\frac{K_0-\e(w')^2}{w^2}(\bargw(X,X)\bargw(Y,Y)-\bargw(X,Y)^2)
$$
for vector fields $X$, $Y$ tangent to $M$.  

The special cases we need are summarized in the following table.  The
first column gives the isometry class of the factor $M$, the next two
columns give the warping function $w$ where the parameter $\ga$ is
defined by the equation in the second column.  The metric $\bargw$ is
then given by equation~(\ref{eq:warp-metric}) with $\e$ as in the
fourth column.  The isometry class of the warped product $\am$ is then
given in the fifth column.  The last two columns give the isometry
class of the hypersurface $M_{t}:=M\times\{t\}$ and the Weingarten
map $S_t$ of $M_t$ with respect to the normal $\eta:=\partial_t$.

\begin{table}[htp]
\caption[]{}
{\footnotesize
\begin{center}
\begin{tabular}{|c|c|c|c|c|c|c|}\hline
$M$&$\ga$&$w$&$\e$&\vvvv$\am$&$M_t$&$S_t$\\
\hline\hline
\vvvv$\r^{n-1}_0(0)$&---&$e^{-t}$&$+1$&$\r^{n}_0(-1)$&$\r^n_0(0)$&$I$\\
\hline
\vvvd$\r^{n-1}_0(0)$&---&$e^{t}
	$&$-1$&\tb{\vvv Part of}{$\r^{n}_1(+1)$}&$\r^{n-1}_0(0)$&$-I$\\
\hline
\tb{$\r^{n-1}_0(K_0)$,}{\vvv$K_0\ge1$}
        &\tb{$\cos(\ga)$}{$=\dfrac{1}{\sqrt{K_0}}$}&
	$\dfrac{\cos(t+\ga)}{\cos(\ga)}$&$+1$&$\r^{n}_0(+1)$
	&$\r^{n-1}_0\!\!\(\dfrac{+1}{\cos^2(t+\ga)}\)$&$\tan(t+\ga)I$\\
\hline
\tb{$\r^{n-1}_0(K_0)$,}{\vvv$-1\le K_0<0$}
	&\tb{$\cosh(\ga)$}{$=\dfrac{1}{\sqrt{|K_0|}}$}
	&$\dfrac{\cosh(t+\ga)}{\cosh(\ga)}$&$+1$&$\r^{n}_0(-1)$
	&$\r^{n-1}_0\!\!\(\dfrac{-1}{\cosh^2(t+\ga)}\)$&$-\tanh(t+\ga)I$\\
\hline
\tb{$\r^{n-1}_0(K_0)$,}{\vvv$0< K_0\le1$}
	&\tb{$\cosh(\ga)$}{$=\dfrac{1}{\sqrt{K_0}}$}
	&$\dfrac{\cosh(t+\ga)}{\cosh(\ga)}$&$-1$&$\r^{n}_1(+1)$
	&$\r^{n-1}_0\!\!\(\dfrac{1}{\cosh^2(t+\ga)}\)$&$-\tanh(t+\ga)I$\\
\hline
\tb{$\r^{n-1}_0(K_0)$,}{\vvv$K_0\le-1$}
        &\tb{$\cos(\ga)$}{$=\dfrac{1}{\sqrt{|K_0|}}$}
    &$\dfrac{\cos(t+\ga)}{\cos(\ga)}$&$-1$&\tb{\vvv Part of}{$\r^{n}_1(-1)$}
	&$\r^{n-1}_0\!\!\(\dfrac{-1}{\cos^2(t+\ga)}\)$&$\tan(t+\ga)I$\\
\hline
\end{tabular}
\end{center}}
\end{table}

Two entries in this table require some explanation.  The first is the
second row which is the warped product of the Euclidean space $\r^{n-1}$
with its standard metric $g=|dx|^2$ and $w(t)=e^{-t}$ so that
$\bargw=e^{-2t}|dx|^2 -dt^2$. If we do the change of variable $y=e^{-t}$,
then the metric becomes
$$
\bargw=\frac{|dx|^2-dy^2}{y^2}
$$
defined on the upper half space $\{(x,y)\in \r^{n-1}\times \r:y>0\}$
(where $|dx|^2:=dx_1^2+\cdots+dx_{n-1}^2$ is the standard metric 
on $\r^{n-1}$).
This is diffeomorphic to $\r^n$ and thus can not be all of the
de~Sitter space $\r^n_1(+1)$ as de~Sitter space is diffeomorphic to
$S^{n-1}\times \r$. Therefore the metric $\bargw$ is not geodesically
complete.  (When $n=2$ this can be seen more concretely by noting that
$c(t)=(\sinh(t),\cosh(t))$ is a geodesic in the upper half plane with
the metric $\bargw$, but that the length of $c$ is
$\int_{-\infty}^\infty\frac{dt}{\cosh(t)}=\pi<\infty$.)  We also
note that if $a\in \r^{n-1}$ and $c>0$ then the metric
$\bargw$ is invariant under the transformations $(x,y)\mapsto
(x+a,y)$ and $(x,y)\mapsto (cx,cy)$ and thus the upper half plane has
a transitive group of isometries.  Thus unlike the Riemannian case a
semi-Riemannian homogeneous space need not be geodesically complete.   
It is convenient to give this half space part of the de~Sitter space a 
name. 
\begin{definition}\sl \label{def:half-space}
Denote by $(\h^n_1(+1),\hm)$ the space $\h^n_1(+1)=\r^{n-1}\times\r$ with
the metric $\hm:=e^{2t}g^{n-1}_0-dt^2$.~\done  
\end{definition}
The other entry that requires some comment is the last row. To simplify
things let $\alpha=0$ and $K_0=-1$ so that $\bargw=\cos^2(t)g-dt^2$ 
where $g$ is the standard metric on the Riemannian manifold
$\r^{n-1}_0(-1)$ of constant sectional curvature~$-1$.  This metric
has coordinate singularities at $t=\pm \pi/2$ but it is
isometric to an open strip in the geodesically complete 
anti--de~Sitter space $\r^n_1(-1)$.
It is also convenient to give this space
a name:
\begin{definition}\sl\label{def:strip}
Denote by $(\els^n_1(-1),\elm)$ the space 
$\els^n_1(-1):=\r^{n-1}_0(-1)\times
(-\pi/2,\pi/2)$ with the metric $\elm:=\cos^2(t)g-dt^2$ where $g$ is the
standard metric on $\r^{n-1}_0(-1)$.~\done  
\end{definition}

Let $\gamma:\r\to \r^n_1(-1)$ be a unit speed timelike geodesic in the
anti--de~Sitter space.  For each $t\in\r$ let $M_t$ be the totally
geodesic hypersurface of $\r^n_1(-1)$ through $\gamma(t)$ and
orthogonal to $\gamma'(t)$.  The set $\{M_t:t\in \r\}$ forms a smooth
foliation of $\r^n_1(-1)$.  This lets us define the time function of
the geodesic $\gamma$ by
\begin{definition}\sl\label{def:time-fcn}
Let $\gamma:\r\to \r^n_1(-1)$ be a unit speed timelike geodesic.  Then
the {\em time function\/} defined by $\gamma$ is the function
$\tau:\r^n_1(-1)$ given by 
$$
\tau(x)=t\quad\mbox{if $x\in M_t$}.
$$
\vskip-22pt~\done
\end{definition}

\subsection{Two Dimensional Rigidity Results}\label{subsec:two-dim}

\newcommand{\un}{\bold{n}}
\newcommand{\rpt}{\bold{RP}^2}

The results in this section are lemmas that are needed to prove our higher
dimensional rigidity results.  

Let $(M^2,g)$ be a semi-Riemannian surface. Let $\nabla$ denote the
Riemannian connection on $M^2$. Then the curvature of $(M^2,g)$ is defined
to be $K = \la R(X,Y)Y,X\ra $ where $X,Y$ are orthonormal.
Now assume that $M$ is oriented and let $e_1, e_2$ denote an oriented 
orthonormal frame  defined locally on $M$. 
Let $\sigma^1, \sigma^2$ denote the dual 
one-forms defined in the domain of definition of the moving frame, 
i.e. if $I$ is the identity map on tangent spaces then
$ I = e_1 \sigma^1 + e_2 \sigma^2 $.
Let $\omega_i^j$ denote the connection forms so that 
$\nabla e_i = e_1 \omega^1_i + e_2 \omega^2_i $.
By a standard calculation
$$
K dA = K \sigma^1\wedge \sigma^2 = \eps_1 d \omega^1_2\,,
$$
where $\eps_1 := \la e_1,e_1 \ra = \pm 1$.

\begin{thm}[Two Dimensional Rigidity]\sl\label{thm:twodrigid}
Let $g_0$ denote the standard flat metric on $\r^2_k$ and let $g$ be any
other metric so that 
$g = g_0$ outside of some compact set $C\subset \r^2$ and 
the curvature $K$ of $g$ does not change sign, i.e. $K \geq 0 $
everywhere or $K \leq 0 $ everywhere.
Then $K_g \equiv 0$ and $(\r^2, g)$ is isometric to $(\r^2_k, g_0)$.
\end{thm}

\begin{pf}  Let $e_1^0, e_2^0$ be the standard basis of $\r^2$, i.e. 
$e_1^0 = { \scriptsize\left [ \begin{array}{c} 1 \\ 0 \end{array} \right ]} ,
e_2^0 = {\scriptsize\left [ \begin{array}{c} 0 \\ 1 \end{array} \right ]} $. 
We claim there is a smooth  orthonormal frame  $e_1, e_2$ for $g$ 
so that $e_i = e_i^0$ on the complement of $C$.  If $g$ is positive
definite (i.e.  $k=0$), then $e_1$,~$e_2$ can be constructed by applying 
the Gram-Schmidt orthogonalization process to 
$e_1^0$,~$e_2^0$ with respect to  $g$. 
The same argument works if $g$ is negative definite
($k=2$).  This leaves the case where $k=1$, whence
$g_0(e_1^0,e_1^0)=-1$ and $g_0(e_2^0,e_2^0)=+1$. 

Let $O_g(\r^2)$ be the bundle of oriented $g$~orthonormal frames
tangent to $\r^2$.  That is the fiber at $P\in \r^2$ is the set of
ordered pairs $(v_1,v_2)$ where $g(v_1,v_1)=-1=-g(v_2,v_2)$,
$g(v_1,v_2)=0$, and the orientation of $v_1, v_2$ agrees with the
orientation of $e_1^0,e_2^0$. This fiber is diffeomorphic to two
disjoint copies of $\r$.  As $\r^2$ is contractible the bundle
$O_g(\r^2)$ is trivial (cf.~\cite{Husemoller:fibre-bundles}), that is
equivalent to a product bundle.  A smooth $\r$~valued function defined
on a neighborhood of the closure of $\r^2\setminus C$ extends to
$\r^2$.  Therefore, the section $e_1^0,e_2^0$ of $O_g(\r^2)$ can be
extended from $\r^2\setminus C$ to a section $e_1,e_2$ defined on all
of $\r^2$.

Now let $D$ be a bounded domain with smooth boundary so that $C \subseteq D$. 
Then 
$$
\int_D K dA = \eps_1 
\int_D
d \omega^1_2 = \eps_1 \int_{\partial D} \omega^1_2. 
$$
But $\nabla e_i = \nabla e_i^0 = 0$ on $\r^2 \backslash C $ so
$\omega^1_2 = 0$ on $\r^2 \backslash C$. Thus the last equation
implies $\int_D K dA = 0$.  As $K$ does not change sign this implies
$K \equiv 0$. It is a standard result that $(\r^2,g)$ is isometric to
$(\r^2_k,g_0)$ (cf.~\cite{Wolf:constant}).  We note that it is also
possible to use the form of the Gauss-Bonnet theorem for two
dimensional space times in \cite{Birman-Nomizu:Gauss-Bonnet} to
complete the proof in the case $k=1$.
\end{pf}

The following variant of this will be used in the proof of 
Theorem~\ref{flatrigid}.

\begin{cor}\sl\label{2d_domain}
Let $g_0$ denote the standard flat metric on $\r^2_k$ and let 
$D\subset \r^2_k$ be a connected domain with smooth boundary and compact
closure in $\r^2_k$.  Let $g$ be any other metric on $D$ so that
$g = g_0$ in some neighborhood of $\partial D$ and 
the curvature $K$ of $g$ does not change sign.
Then $K_g \equiv 0$.~\done
\end{cor}

We now turn to rigidity on the sphere.  
\begin{thm}[Toponogov, \cite{top}]\sl Let $(S^2,g)$ be the 
two dimensional sphere
with a metric so that the Gaussian curvature satisfies $1\le K$.  Then any
simple closed geodesic $\gamma$ has length at most $2\pi$, and if the
length of $\gamma$ is $2\pi$ then $(S^2,g)$ is 
isometric to the standard
sphere $(S^2,g_0)$ and $\gamma$ is a great circle.~\done
\end{thm}
We owe both the statement and 
proof of the following result to E. Calabi~\cite{Calabi:Sphere-Rigidity}.
\begin{thm}\sl \label{thm:calabi}
Let $(S^2,g)$ be the two dimensional sphere with
a metric of class $C^{1,1}$ whose Gaussian curvature satisfies 
$ 0\le K\le 1$.  
Then any simple closed geodesic $\gamma$ on $(S^2,g)$ has length at least 
$2\pi$. If the length of $\gamma$ is $2\pi$, then either $(S^2,g)$ is isometric
to the standard round sphere $(S^2,g_0)$ and $\gamma$ is a great circle on
$(S^2,g_0)$ or $(S^2,g)$ is isometric to a circular cylinder of
circumference $2\pi$ capped by two unit hemispheres and $\gamma$ is a belt
around the cylinder. Thus if $K$ is continuous (for example
when $g$ is at least $C^2$) or if $K>0$ then $(S^2,g)$ is isometric to
the standard round sphere.
\end{thm}
The lower bound on the length of a simple closed closed
geodesic is well known (cf.~Remark~\ref{rk.len}), it is the rigidity result
that is of interest here.

\begin{lemma}\sl Let $k(t)$ be an $L^\infty$ function on $[0,\infty)$ so that 
$0\le k(t)\le 1$ and let $y(t)$ be defined by the initial value problem 
$$
y''(t)+k(t)y(t)=0,\quad y(0)=1,\quad y'(0)=0.
$$
Denote the smallest positive zero of $y$ by $\beta$ (it may be that
$\beta=\infty$).  Then 
$0\le -y'(t)\le 1$ for $0\le t\le \beta$.  If $y'(t_0)=-1$ for some 
$t_0\in [0,\beta]$, then $t_0=\beta<\infty$, $\beta\ge\pi/2$ and 
$$
y(t)=\left\{\begin{array}{rl}1,
& 0\le t < \beta-\pi/2\\ \cos(t-(\beta-\pi/2)),& \beta-\pi/2<t\le\beta
\end{array}\right.,\quad
k(t)=\left\{\begin{array}{rl}0,
& 0\le t < \beta-\pi/2\\ 1,& \beta-\pi/2<t\le\beta.
\end{array}\right.
$$
If $k$ is continuous, for example if $y$ is $C^2$, and $y'(t_0)=-1$ with
$t_0\le \beta$, then $t_0=\beta=\pi/2$, $k\equiv1$ and $y(t)=\cos(t)$ on
$[0,\pi/2]$.
\end{lemma}

\begin{pf} Note  on the interval $[0,\beta)$ that $(y')'=y''=-ky\le0 $
as $y>0$ and $k\ge 0$ on $[0,\beta)$. Thus $y'$ is monotone decreasing on 
$[0,\beta]$.  As $y'(0)=0$ this implies $y'\le0$ on $[0,\beta]$. Thus 
\begin{equation}
\left(y^2+(y')^2\right)'=2yy'+2y'y''=2yy'-2y'ky=2yy'(1-k)\le 0 \label{decr}
\end{equation}
on $[0,\beta)$ as $k\le 1$ and $y'\le 0$.  Using the initial conditions 
for $y$ and continuity the last inequality implies
\begin{equation}
y^2+(y')^2\le1\quad \mbox{on} \quad [0,\beta].  \label{sumsq}
\end{equation}
These inequalities imply $0\le -y'\le 1$ on $[0,\beta]$. 

If $t_0\in[0,\beta]$ and $y'(t_0)=-1$, then the inequality~(\ref{sumsq})
implies $y(t_0)=0$.  But from the definition of $\beta$ as the smallest
positive zero of $y$ this implies $t_0=\beta$. Then
$y(\beta)^2+y'(\beta)^2=1$ and equation (\ref{decr}) yields
$y'(1-k)\equiv0$ on $[0,\beta)$.  As $y'$ is monotone decreasing on 
$[0,\beta)$, there is a point $t_1\in[0,\beta)$ so that 
$y'\equiv0$ on $[0,t_1]$ and $0>y'>-1$ on $(t_1,\beta)$.  Then on 
$[0,t_1]$ we have $y\equiv1$ and $k\equiv0$ (as $ky=-y''=0$). Also 
$y'(1-k)\equiv 0 $ and $y'\ne0$ on $(t_1,\beta)$ implies $k\equiv1$ on 
$(t_1,\beta)$.  But $y(t_1)=1$ and $y'(t_1)=0$ so $y(t)=\cos(t-t_1)$ 
on $(t_1,\beta)$.  As $y(\beta)=0$ this implies $t_1=\beta-\pi/2$ 
on $(t_1,\beta)$.  This completes the proof.
\end{pf}

\begin{pf*}{Proof of the theorem} Let $c:[0,L]\to S^2$ be a unit
speed parameterization of the closed geodesic $\gamma$, and let $\un$ be a
unit normal along $c$.  For each $s\in [0,L]$ let $\beta(s)$ be the cut
distance from the curve $\gamma$ along the geodesic $t\mapsto
\exp_{c(s)}(t\un(s))$.  Define map $F(s,t)$ on the set of ordered pairs 
$(s,t)$ with $s\in[0,L]$ and $0\le t\le\beta(s)$ by 
$$
F(s,t)=\exp_{c(s)}(t\un(s)),\qquad 0\le s\le L, \quad 0\le t\le\beta(s).
$$
Then $s,t$ are Fermi coordinates on the disk $M$ bounded by $\gamma$ and
with inner normal $\un$.  In these coordinates the metric $g$, Gaussian curvature
$K$ and the area form $dA$ are given by
$$
g=E^2\,ds^2+ dt^2,\quad K=\frac{- E_{tt}}{E},\quad dA=E\,ds\,dt.
$$
And because $c$ is a geodesic $E(s,0)\equiv1$ and $E_t(s,0)\equiv0$.
Thus for fixed $s$ the function $y(t):=E(s,t)$ satisfies 
$y''+Ky=0$, $y(0)=1$, and $y'(0)=0$ as in the lemma.

Now apply the Gauss-Bonnet theorem to the disk $M$.  As the boundary is a
geodesic, the boundary term of the formula drops out:
\begin{align*}
2\pi&=\int_MK\,dA
	=\int_0^L\int_0^{\beta(s)}- E_{tt}\,dtds\\
	&=\int_0^L(-E_t(s,\beta(s)))\,ds\qquad\quad\mbox{(as $E_t(s,0)=0$)}\\
	&\le\int_0^L1\,ds \qquad\mbox{\ \hskip1in (by the lemma)}\\
	&=L.
\end{align*}
which proves the required lower bound on the length of $\gamma$.  If
$L=2\pi$, then $E_t(s,\beta(s))=-1$ for all $s\in[0,L]$. Again by the lemma
in the coordinates $s,t$ on $M$
$$
K(s,t)=\left\{\begin{array}{rl}0,&\quad 0\le t<\beta(s)-\pi/2,\\
	1,&\quad \beta(s)-\pi/2<t\leq  \beta(s).
	\end{array}\right.
$$
Let $M_{+1}$ denote the interior of the set $\{x\in M:K(x)=+1\}$ so that 
$M_{+1}=\{\exp_{c(s)}(t\un(s)):s\in[0,2\pi],\beta(s)-\pi/2<t\le\beta(s)\}$. 
Let $s_0\in[0,2\pi]$ be a point where $\beta(s)$ is maximal.  Then the open 
disk $B(x_0,\pi/2)$ of radius $\pi/2$ about 
$x_0:=\exp_{c(s_0)}{\beta(s_0)\un(s_0)}$ is
contained in $M_{+1}$, for if not it would meet $\partial M_{+1}$ at some
point $\exp_{c(s)}((\beta(s)-\pi/2)\un(s))$ and this point is a distance of 
$\beta(s)-\pi/2$ from $\gamma$.  Thus the distance of
$x_0=\exp_{c(s_0)}(\beta(s_0)\un(s_0))$ to $\gamma$ is less than
$\pi/2+(\beta(s)-\pi/2)=\beta(s)$, which contradicts the maximality of
$\beta(s_0)$. Thus $B(x_0,\pi/2)\subseteq M_{+1}$.  But using the
Gauss-Bonnet theorem and  $K\equiv+1$ on $M_{+1}$
$$
2\pi\ge\int_{M_{+1}}K\,dA=\mathop{\mathrm{Area}}(M_{+1})\ge
	\mathop{\mathrm{Area}}(B(x_0,\pi/2))=2\pi.
$$
So $M_{+1}=B(x_0,\pi/2)$.  From this it follows
$s\mapsto\beta(s)$ is constant and thus the disk $M$ bounded by $\gamma$
and with inner normal $\un$ is a cylinder of circumference $2\pi$ capped
at one end with a hemisphere.
The same argument applied to the disk bounded by $\gamma$ and having $-\un$
as inward normal shows  $(S^2,g)$ is two of these capped cylinders glued
together along $\gamma$, which is equivalent to the statement of the
theorem.
\end{pf*}

\begin{remark}\normalshape\label{rk.len} If one is only interested in the length
of closed geodesics, there are higher dimensional versions of 
Theorem~\ref{thm:calabi}. If
$(M,g)$ is a compact orientable Riemannian manifold of even dimension with
sectional curvatures satisfying, $0<K_M\le1$, then Klingenberg has shown
every closed geodesic has length $\ge2\pi$.  This is equivalent to his
well known lower bound on the injectivity radius of compact oriented even
dimensional manifolds.  For a proof see the 
book~\cite[Chapter 5]{Cheeger-Ebin}.
For odd dimensional manifolds another theorem of Klingenberg's implies if 
$(M,g)$ is a compact simply connected manifold of whose sectional
curvature satisfies $1/4<K_M\le 1$, then any closed geodesic of $(M,g)$
has length at least $2\pi$.  Again a proof can be found in   
\cite[Chapter 5]{Cheeger-Ebin}. (The original proofs of Klingenberg are in 
\cite{Klingenberg:Contributions,Klingenberg:Conjugate}.) We also note that in dimension~3 for any $\varepsilon>0$
there are examples of metrics $g$ on $M=S^3$ (due to Berger) so that the
sectional curvatures satisfy $1/9-\varepsilon \le K_M\le 1$, but 
$(S^3,g)=(M,g)$ has a geodesic of length less than $2\pi$.
Cf.~\cite[Example 3.35, page~70]{Cheeger-Ebin}.  
To the best of our knowledge there
is no known rigidity result in these theorems.~\done
\end{remark}

\section{Basic Comparison Results}\label{sec:basic}

Most  of the results in this section will be stated as results about
systems of ordinary differential equations (usually in the form of
solutions to matrix valued differential equations).  These results
then apply to the second fundamental forms (or Weingarten maps) of
parallel families of hypersurfaces in a Semi-Riemannian manifold.     

\subsection{Comparisons for Parallel Hypersurfaces}

In this section we will prove a comparison theorem for solutions of the 
matrix Riccati equation.

Let $E$ be an $n$~dimensional real inner product space with an inner
product of index~$k$. The inner product will be denoted by $\la\,,\ra$.
A linear map $A:E\to E$ is {\em self-adjoint\/} iff 
$\la Ax,y\ra=\la x,Ay\ra$  for all $x,y\in E$.  Unlike the
case of positive definite spaces a self-adjoint linear map need not have 
real eigenvalues.  (An easy example of this is $E=\r^2$,
$\la\,,\ra=dx^2-dy^2$
and $A={\scriptsize\left[\begin{array}{cc}0&1\\-1&0\end{array}\right]}$. 
Then the eigenvalues of $A$ are $\pm\sqrt{-1}$.)

A self-adjoint linear map $A$ is {\it positive definite\/} if and only if 
$\la Ax,x\ra>0$ for all $x\ne 0$.  The eigenvalues of a positive definite
linear map are real. (This can be seen by using the principal axis theorem
to diagonalize the inner product $\la\,,\ra$ with respect to the positive
definite inner product $\la A\,, \ra$.)  However the eigenvalues of
a positive definite map need not be positive.  In fact if the index of
the inner product is $k$ then a positive definite map will have exactly
$k$ negative eigenvalues, however having $k$ negative eigenvalues is not
enough to insure   $A$ is positive definite.  Unlike 
the positive definite case the identity map is not positive definite. A 
linear map is {\em positive semi-definite\/} iff $\la Ax,x\ra\ge 0$ for 
all $x\in E$.  If $A$ and $B$ are self-adjoint, then we write
$A\le B$ iff $B-A$ is positive semi-definite, and $A<B$ iff $B-A$ is
positive definite. The set of positive definite maps forms a cone.

\begin{lemma}\sl \label{kernel}
If $A\ge 0$ and $\la Ax_0,x_0\ra=0$, then $Ax_0=0$.
\end{lemma}

\begin{pf} Because $A$ is positive semi-definite and $\la A x_0,x_0\ra=0$ 
we have for all $h\in E$ that 
$0\le \la A(x_0+h),(x_0+h)\ra=2\la Ax_0,h\ra+\la Ah,h\ra$. 
This can only be non-negative for all $h$ if $\la Ax_0,h\ra=0$ for all 
$h$.  This implies  $Ax_0=0$.
\end{pf}

For $i=1,2$ let $t\mapsto R_i(t)$ be smooth maps from the real numbers
$\r$ to the space of self-adjoint linear maps on $E$.  Also let 
$A_1$ and $A_2$ be self-adjoint linear maps on $E$ and denote by $I$
the identity operator on $E$.  Define two maps $t\mapsto F_i(t)$ from
the reals to the space of linear maps on $E$ by the initial value problems,
\begin{equation}
F_i''(t)+R_i(t)F_i(t)=0,\quad F_i(0)=I,\quad F_i'(0)=-A_i\,; \qquad i=1,2 .
\label{jacobi}
\end{equation}
At the points $t$ where $F_i(t)$ is invertible set
$$
S_i(t)=-F_i'(t)F_i(t)^{-1}\,; \qquad i=1,2.
$$
Let $F_i(t)^*$ be the transpose of $F_i(t)$ with respect to the inner
product $\la\,,\ra$ (defined by $\la F_i(t)^*x,y\ra=\la x,F_i(t)y\ra$).
Then using the equation (\ref{jacobi}) it follows 
$$
\frac{d}{dt}(F_i(t)^*F_i'(t)-F_i'(t)^*F(t))=0\,;\qquad i=1,2.
$$
which implies   the linear maps $S_i(t)$ are self-adjoint at all points
where they are defined. It also follows directly from (\ref{jacobi}) that
in some neighborhood of $0$, the $S_i(t)$ satisfy the initial value
problem
$$
S_i'(t)=S_i(t)^2+R_i(t),\quad S_i(0)=A_i\,; \qquad i=1,2.
$$

\begin{thm}\sl \label{hypercomp} With the notation above, if 
both $S_1(t)$ and $S_2(t)$ are defined on all of the interval $[0,b]$ and 
if
$$
A_1\le A_2,\qquad R_1(t)\le R_2(t)\  \mbox{for all}\  t\in [0,b],
$$
then $S_1(t)\le S_2(t)$ on $[0,b]$.  If $S_1(b)=S_2(b)$, then
$A_1=A_2$ and $R_1(t)\equiv R_2(t)$ on $[0,b]$.
\end{thm}

\begin{pf}  We first prove this under the assumption  $A_1<A_2$ and 
$R_1(t)<R_2(t)$ for all $t\in [0,b]$ and show   in this case  
$S_1(t)<S_2(t)$ for $t\in [0,b]$.  Assume, toward a contradiction,  
this is false.  Then there is a smallest $t_0$ so that $S_2(t_0)-S_1(t_0)$
is not positive definite.  By continuity $S_2(t_0)-S_1(t_0)$ is positive
semi-definite.  As it is not positive definite there is a nonzero vector
$x_0\in E$ so that $\la (S_2(t_0)-S_1(t_0))x_0,x_0\ra=0$.  By Lemma
\ref{kernel} this yields  $S_1(t_0)x_0=S_2(t_0)x_0$.

Let $f(t):=\la (S_2(t)-S_1(t))x_0,x_0\ra$.  Then, using 
$S_1(t_0)x_0=S_2(t_0)x_0$,
\begin{align*}
f'(t_0)&=\la (S_2'(t_0)-S_1'(t_0))x_0,x_0\ra\\
	&=\la S_2(t_0)x_0,S_2(t_0)x_0\ra-\la S_1(t_0)x_0,S_1(t_0)x_0\ra
	  +\la (R_2(t_0)-R_1(t_0))x_0,x_0\ra\\
	&=\la (R_2(t_0)-R_1(t_0))x_0,x_0\ra\\
	&>0.
\end{align*}
But $f$ is positive on $[0,t_0)$ and $f(t_0)=0$ thus $f'(t_0)\le 0$.  This
contradicts the last equation and completes the proof in the case 
$A_2-A_1$ and $R_2(t)-R_1(t)$ are positive definite.

For the general case let $B$ be any positive definite matrix.  
Then for small $\delta >0$ define $S_\delta$ by the initial value problem
$$
S_\delta'(t)=S_\delta(t)^2 +R_2(t)+\delta B,\quad S_\delta(0)=A_2+\delta B.
$$
By what has already been done we have $S_1(t)< S_\delta(t)$ on $[0,b]$.
Using that $S_\delta $ depends continuously on $\delta$,
$$
S_2(t)=\lim_{\delta\downarrow 0} S_\delta(t) \ge S_1(t).
$$
This proves the inequality in the general case. 

Assume   $S_1(b)=S_2(b)$ and let $t_0\in [0,b)$.  Then define $S$ 
on $[t_0,b]$ by 
$$
S'(t)=S(t)^2+R_1(t),\quad S(t_0)=S_2(t_0)
$$
so $S$ satisfies the same differential equation as $S_1$ and agrees
with $S_2$ at $t_0$.  By what has already been shown this implies 
$S_1(t)\le S(t)\le S_2(t)$ on $[t_0,b]$.  As $S_1(b)=S_2(b)$
this yields $S(b)=S_1(b)$.  But $S_1$ and $S$ satisfy the same 
first order equation and agree at the point $b$, so $S\equiv S_1$ on 
$[t_0,b]$.  In particular 
$S_1(t_0)=S(t_0)=S_2(t_0)$.  As $t_0$ was arbitrary $S_1\equiv S_2$.
By uniqueness this
implies $A_1=A_2$ and $R_1\equiv R_2$.
\end{pf}

A standard part of the comparison theory in the
positive definite case is, with the notation above,  if $ S_2(t)$ 
is defined on all of $[0,b]$ then so is $S_1(t)$.  This does not hold on 
in the general indefinite metric case.  We give two examples on $\r^2$ with
the inner product $\la\,,\ra=dx^2-dy^2$.  First let $A_1=A_2=0$, 
$R_1(t)\equiv 0$ and 
$R_2(t)\equiv {\scriptsize\left[\begin{array}{cc}1&0\\0&0\end{array}\right]}$.  
Then $R_1\le R_2$, $S_1=0$ is defined on all of $[0,\infty)$ but 
$S_2(t)={\scriptsize\left[\begin{array}{cc}\tan
t&0\\0&0\end{array}\right]}$ 
is only defined on $[0,\pi/2)$. In the second example again let $A_1=A_2=0$, 
but let 
$R_1(t)\equiv
{\scriptsize\left[\begin{array}{cc}0&0\\0&1\end{array}\right]}$ 
and $R_2(t)\equiv 0$. Again $R_1(t)\le R_2(t)$. 
$S_1(t)={\scriptsize\left[\begin{array}{cc}0&0\\0&\tan
t\end{array}\right]}$ 
is only defined on $[0,\pi/2)$ but $S_2(t)=0$ is defined on all of $[0,\infty)$.

However if there is a two sided curvature bound, then there is a
an estimate of the size of the domain of definition.

\begin{prop}\sl \label{twosided} Let $R_1(t), R_2(t), R_3(t)$ be continuous 
functions of $t\in \r$ with values in the symmetric linear maps on 
$(E,\la\,,\ra)$ so that $R_1\le R_2\le R_3$ on $[0,b]$. For $i=1,2,3$,
define $S_i(t)$ 
by the initial value problems $S_i'=S_i^2+R_i$ and assume  
$S_1(0)\le S_2(0)\le S_3(0)$. If $S_1$ and $S_3$ are defined on all of 
$[0,b)$, then so is $S_2$.
\end{prop}

\begin{pf} Let $[0,t_0)$ be the maximal interval of definition of $S_2$.  If
$t_0 \geq b$  we are done, so assume, toward a contradiction,  
$t_0 < b$.  
By the comparison theorem $S_1(t)\le S_2(t)\le S_3(t)$ on $[0,t_0)$.  As 
$S_1$ and $S_3$ are continuous on $[0,t_0]$ their ranges are bounded.  Thus
there are self-adjoint $A$ and $B$ so that $A\le S_1(t)\le S_3(t)\le B$.  But
the set ${\cal C}=\{C: A\le C\le B\}$ is a compact
subset of the space of 
self-adjoint linear maps and $S_2(t)\in {\cal C}$ for $t\in [0,t_0)$.  
Therefore there is a sequence $t_k\nearrow t_0$ so that
$\lim_{k\to\infty}S_2(t_k)$
exists, say $\lim_{k\to\infty}S_2(t_k)=G$. Choose any norm $\|\cdot\|$ on 
the space of self-adjoint linear maps on $E$.   Then, as $S_2(t)$ stays in 
a compact set, there is a constant $c_0$ so that $\|S_2(t)^2+R_2(t)\|\le c_0$ 
for all $t\in [0,t_0)$.  Thus the differential equation for $S_2$ implies 
$\|S_2(t)-S_2(s)\|\le c_0|s-t|$ for all $t,s\in [0,t_0)$.  Then 
$\lim_{t\nearrow t_0}=\lim_{k\to \infty}S_2(t_k)=G$. Define $S(t)$ by 
$S'=S^2+R_2$ and $S(t_0)=G$.  Then by uniqueness of solutions to initial
value problems we find that 
$S_2=S$ in a neighborhood of $t_0$, contradicting that 
$[0,t_0)$ was the maximal interval of definition of $S_2$
\end{pf}

\subsubsection{Comparisons for Curvature Tensors}

For the later applications we need not just the above comparison theorems
for the Weingarten map, but also comparison theorems
for the intrinsic  curvature tensor.
The idea is to use the Gauss curvature equation to relate the estimates on
Weingarten map to the sectional curvature.  We start by giving some linear
algebraic results.

Recall  the inner product $\la\,,\ra$ induces an inner product 
(also denoted by $\la\,,\ra$) on $\bigwedge^2E$.  On decomposable
elements it is given by,
$$
\la x_1\w x_2,y_1\w y_2\ra
	=\la x_1,y_1\ra\la x_2,y_2\ra-\la x_1,y_2\ra\la x_2,y_1\ra.
$$
Also a linear map $A:E\to E$ induces a linear map 
$\W^2(A):\W^2E\to\W^2E$ given on decomposable elements by 
$\W^2(A)x\w y=Ax\w Ay$.

\begin{definition}\sl If $A$, $B$ are self-adjoint linear maps on $E$ then 
define $\W^2(A)\preceq \W^2(B)$ to mean that 
$\la \W^2(A)x\w y,x\w y\ra\le \la \W^2(B)x\w y,x\w y\ra$
for all decomposable $x\w y$.  This is equivalent to 
$$
\la Ax,x\ra\la Ay,y\ra-\la Ax,y\ra^2\le \la Bx,x\ra\la By,y\ra-\la Bx,y\ra^2
$$
holding for all $x,y\in E$.~\done
\end{definition}

\begin{lemma}\sl\label{positive} 
If $A$ is self-adjoint and $\la\W^2(A)x\w y,x\w y\ra=0$ for all
decomposable vectors $x\w y$, then $A$ has rank one and so $\W^2(A)=0$.
If $A$ is self-adjoint and $0\le A$ (or $A\le 0$) then
$\W^2(A)\succeq 0$.
\end{lemma}

\begin{pf} If $\la\W^2(A)x\w y,x\w y\ra=0$ then 
$$
\la Ax,x\ra\la Ay,y\ra-\la Ax,y\ra^2=0.
$$
Let $h\in E$.  Then replace $x$ by $x+th$ in the last equation and take 
$\frac{d}{dt}|_{t=0}$ to get  
$$
2\la Ax,h\ra\la Ay,y\ra-2\la Ax,y\ra\la h,Ay\ra=0
$$
for all $x,y,h\in E$.  If $A=0$ there in nothing to prove.  Thus assume  
$A\ne 0$.  Then, as $A$ is self-adjoint, there is a $y\in E$ so that 
$\la Ay,y\ra\ne 0$.  From the last equation we see  if $\la h,Ay\ra=0$,
then also $\la Ax,h\ra=0$.  This implies  $Ax$ is a scalar multiple of 
$Ay$.  As $x$ was arbitrary, this shows $Ay$ spans the range of $A$,
so  $A$ has rank one.

If $A\ge0$ (or $A \leq 0$), then the Cauchy-Schwartz inequality applied to 
the inner product $\la A\cdot,\cdot\ra$ implies $\la Ax,x\ra\la Ay,y\ra-\la
Ax,y\ra^2\ge0$, i.e. $\W^2(A)\succeq0$.
\end{pf}

\begin{lemma} \label{defcomp}
\sl If $A$, $B$ are self-adjoint and  $0\le A\le B$
(or $B\le A\le 0$), then $\W^2(A)\preceq \W^2(B)$. If also either $A$ or
$B$ is positive definite and 
$\la\W^2(A)x\w y,x\w y\ra=\la\W^2(B)x\w y,x\w y\ra$ for all $x, y\in E$
then $A=B$.
\end{lemma}

\begin{pf} We first assume  $B$ is positive definite.  Let $x,y\in E$ be
linearly independent and let $V=\mathop{\mathrm{Span}}(x,y)$.  Let $\alpha$ and 
$\beta$ be the inner products defined on $V$ by $\alpha(u,v)=\la Au,v\ra$,
$\beta(u,v)=\la Bu,v\ra$.  As $B$ is positive definite so is $\beta$.  
Thus by the principal axis theorem the forms $\alpha$ and $\beta$ can be 
simultaneously diagonalized.  Therefore there is a basis $x_1$, $y_1$ of 
$V$ with $x_1\w y_1=x\w y$ and $\la Ax_1,y_1\ra=\la Bx_1,y_1\ra=0$.
But for all $u\in E$, $0\le\la Au,u\ra\le \la Bu,u\ra$ so  
\begin{align*}
\la {\W}^2(A)x\w y,x\w y\ra&=\la Ax_1,x_1\ra\la Ay_1,y_1\ra \\
&\leq \la Bx_1,x_1\ra\la By_1,y_1\ra=\la {\W}^2(B)x\w y,x\w y\ra .
\end{align*}
If equality holds then $\la Ax_1,x_1\ra=\la Bx_1,x_1\ra$ and 
$\la Ay_1,y_1\ra=\la By_1,y_1\ra$. As $x_1, y_1$ diagonalizes both
$\alpha$ and $\beta$ this implies  $\la Au,u\ra=\la Bu,u\ra$ for all
$u\in V$.  So if $\la\W^2(A)x\w y,x\w y\ra=\la\W^2(B)x\w y,x\w y\ra$ for 
all $x, y\in E$, then $\la Au,u\ra=\la Bu,u\ra$ for all $u\in E$, which
implies $A=B$.  

A similar argument works if $A$ is positive definite.  Finally if $B$ is only 
positive semi-definite then let $C$ be positive definite and let 
$B_\ell=B+\frac1\ell C$.  Then $B_\ell$ is positive definite and
$\lim_{\ell\to\infty}B_\ell =B$.  But $0\le A\le B_\ell$ implies 
$\W^2(A)\le \W^2(B_\ell)$ as $B_\ell$ is positive definite, so a limit 
argument shows this also holds when $B = \lim_{\ell \to \infty} B_\ell$ 
is positive semi-definite.
\end{pf}

\begin{lemma}\sl \label{rank1}
A rank one self-adjoint linear map $S:E\to E$ 
is of one of the two forms $Sx=\la x,e\ra e$ or $Sx=-\la x,e\ra e$, 
for some $e \in E$.
\end{lemma}
\begin{pf} A rank one linear map is of the form $Sx=\la x,u\ra v$ with both 
$u$, $v$ non-zero.  The map $S$ is self-adjoint iff 
$\la Sx,y\ra=\la x,Sy\ra$, which in our case this reduces to 
$\la x,u \ra \la y,v \ra=\la x,v\ra\la y,u\ra$.  This implies  
$u$ and $v$ are linearly dependent.  As they are both non-zero 
$v=\lambda u$ for some non-zero real number $\lambda$.  If $\lambda>0$
let $e=\sqrt{\lambda}\,u$.  Then $Sx=\la x,u\ra v=\la x, e\ra e$.
If $\lambda <0$ set $e=-\sqrt{|\lambda|}\,u$.  Then 
$Sx=\la x,u\ra v=-\la x, e\ra e$.
\end{pf}

\begin{thm}\sl \label{curvcomp}
Let $t\mapsto S(t)$ be a smooth map from $[0,b]$ to the 
self-adjoint maps on $E$.  Assume  $S$ satisfies
$$
S'(t)=S(t)^2+R(t),\quad S(0)=0 ,
$$
where $R(t)\ge 0$ on $[0,b]$ (or $R(t)\le 0$ on $[0,b]$).  Then 
$\W^2(S(t))\succeq0$ on $[0,b]$.  If $\la \W^2(S(b))x\w y,x\w y\ra=0$
for all $x$, $y$ then there is a self-adjoint rank one map $P$ and smooth 
functions $u,r:[0,r]\to \r$ so that 
$$
S(t)=u(t)P,\quad R(t)=r(t)P,\quad \mbox{\sl and} \quad u'(t)=u^2(t)+r(t).
$$
\end{thm}

\begin{pf} We deal with the case  $R(t)\ge 0$ on $[0,b]$, the case with
$R(t)\le 0$ being similar.  By the comparison Theorem~\ref{hypercomp} 
$S(t)\ge 0$ on $[0,b]$.  By Lemma~\ref{positive} this implies 
$\W^2(S(t))\succeq0$.  Now assume  $\la \W^2(S(b))x\w y,x\w y\ra=0$
for all $x$, $y$.  Then again by Lemma~\ref{positive}, $S(b)$ has rank one or
less.  If $S(b)=0$ then Theorem~\ref{hypercomp} implies $S(t)\equiv 0$ and
$R(t)\equiv 0$. 

Now assume  $S(b)$ has rank one.  Again by Theorem~\ref{hypercomp}
there is a $t_0\in [0,b)$ so that $S(t)=0$ on $[0,t_0]$ and $S(t)$ has rank
one on $(t_0,b]$.  By Lemma~\ref{rank1} each $t\in (t_0,b]$ there is 
$e(t)\in E$ so that $S(t)=\pm\la\,\cdot\,,e(t)\ra e(t)$.  
Then $S'(t)=\pm(\la\,\cdot\,,e'(t)\ra e(t)+\la\,\cdot\,,e(t)\ra e'(t))$, 
and $S(t)^2=\la\,\cdot\,,e(t)\ra \la e(t),e(t)\ra e(t)$.  From the
differential equation for $S$, 
$$
R(t)=S'(t)-S(t)^2=\pm(\la\,\cdot\,,e'(t)\ra e(t)+\la\,\cdot\,,e(t)\ra e'(t))
-\la\,\cdot\,,e(t)\ra \la e(t),e(t)\ra e(t).
$$
Using  $R(t)\ge 0$, 
$$
\la R(t)x,x\ra=\pm2\la x,e(t)\ra\la x,e'(t)\ra-\la x,e(t)\ra^2\la
e(t),e(t)\ra\ge0
$$
for all $x$.

\noindent
{\bf Claim:}\quad For any $t\in (t_0,b]$ the vectors $e(t)$ and $e'(t)$ are
linearly dependent.

To see this assume  for some $t$ that $e(t)$ and $e'(t)$ are linearly
independent. Then there are vectors $x_0, h\in E$ so that $\la x_0,e(t)\ra=1$,
$\la x_0,e'(t)\ra=0$, $\la h,e(t)\ra=0$, and $\la h,e'(t)\ra=1$.  Let
$\lambda$ be any real number and let $x=x_0+\lambda h$.  Using this $x$ in
the expression above for $\la R(t)x,x\ra$ implies that for all $\lambda$
$$
\la R(t)x,x\ra=\pm 2\lambda -\la e(t),e(t)\ra\ge 0.
$$
which is impossible.  This proves the claim.

As $e(t)$ and $e'(t)$ are linearly dependent for all $t$ the span of 
$e(t)$ stays constant.  Thus there is a constant vector $e_0$ so that
$\mathop{\mathrm{Span}}e(t)=\mathop{\mathrm{Span}}e_0$ for all $t$.  Let $P$ be a
non-zero rank one self-adjoint map with the range of $P$ being the span of
$e_0$.  (Say $P=\la\,\cdot\,,e_0\ra e_0$.) Then $S(t)=u(t)P$ for some
smooth function $u$.  From the equation for $S(t)$ we have
$R(t)=S'(t)-S(t)^2=(u'(t)-u(t)^2)P$ so setting $r(t)=u'(t)-u(t)^2$ completes
the proof.
\end{pf}

\begin{thm}\sl \label{defcompcurv} 
Let $r(t)$ be a smooth function on $[0,b]$ and $a$ a real number.  Assume
there is a function $u$ on $[0,b]$ so that 
$$
u'(t)=u^2(t)+r(t),\quad u(0)=a.
$$
Let $A$ be a positive definite self-adjoint map on $E$, and 
$t\mapsto S(t)$ a map from $[0,b]$ to the self-adjoint maps
on $E$ so that
$$
S'(t)=S(t)^2+R(t),\quad S(0)=aA
$$
and  $R$ satisfies one of the following conditions:
\begin{enumerate}
\item  $R(t)\ge r(t)A$ on $[0,b]$ and $u(b)>0$, or
\item   $R(t)\le r(t)A$ on $[0,b]$ and $u(b)<0$. 
\end{enumerate}
Then $\W^2(S(b))\succeq u(b)^2\W^2(A)$.  If 
$\la \W^2(S(b))x\w y,x\w y\ra=u(b)^2\la \W^2(A) x\w y,x\w y\ra$ 
for all $x,y\in E$ then $S(t)\equiv u(t)A$ and $R(t)\equiv r(t)A$.
\end{thm}

\begin{pf} We prove this under the first assumption, the proof in the second
case being almost identical.  By  Theorem~\ref{hypercomp} we have
$0\le u(b)A\le S(b)$.  By Lemma~\ref{defcomp} this implies 
$\W^2(S(b))\succeq u(b)^2\W^2(A)$.  Also by Lemma~\ref{defcomp} if
$\la \W^2(S(b))x\w y,x\w y\ra=u(b)^2\la \W^2(A) x\w y,x\w y\ra$ 
for all $x,y\in E$ then $S(t)\equiv u(t)A$ and $R(t)\equiv r(t)A$.
Then the uniqueness part of Theorem~\ref{hypercomp} implies 
$S(t)=u(t)A$ and  $R(t)=r(t)A$.
\end{pf}

\begin{thm}\sl \label{defcompcurv2} 
Let $r(t)$ be a smooth function on $[0,b]$ and $a>0$ a real number.  
Assume there is a function $u$ on $[0,b]$ so that 
$$
u'(t)=u^2(t)+r(t),\quad u(0)=a.
$$
Let $A$ be a positive definite self-adjoint map on $E$, and 
$t\mapsto S(t)$ a map from $[0,b]$ to the self-adjoint maps
on $E$ so that
$$
S'(t)=S(t)^2+R(t),\quad S(0)=aA
$$
and assume  
$$
R(t)\le r(t)A,\ \mbox{and}\ S(t)>0\ \mbox{on}\ [0,b].
$$
Then $\W^2(S(b))\le u(b)^2\W^2(A)$.  If 
$\la \W^2(S(b))x\w y,x\w y\ra=u(b)^2\la \W^2(A) x\w y,x\w y\ra$ 
for all $x,y\in E$ then $S(t)\equiv u(t)A$ and $R(t)\equiv r(t)A$.
\end{thm}

\begin{pf} A variant on the proof of the last theorem.
\end{pf}

\begin{remark}\normalshape
In the Riemannian case there are volume comparison theorems for the
volume of geodesic balls and tubes.  In the Lorentzian case there are
volume comparison for hypersurfaces parallel to a spacelike
hypersurface. In the notation of this section these
reduce to giving inequalities between $\det(F_1(t))$ and $\det (F_2(t))$.  
We now give examples to show  there are no such inequalities in the
general indefinite metric case.  We will work on $\r^2$ with the metric 
$\la\,,\ra=dx^2-dy^2$.  
First let $A_1=A_2=0$, $R_1(t)\equiv 0$, and 
$R_2(t)\equiv {\scriptsize \[\begin{array}{cc}1&0\\0&0\end{array}\]}$.
Then $R_1\le R_2$. Define $F_i$ by $F_i''+R_iF_i=0$, $F_i(0)=I$,
$F_i'(0)=A_i=0$.  Then 
$F_1(t)={\scriptsize \[\begin{array}{cc}1&0\\0&1\end{array}\]}$ and
$F_2(t)={\scriptsize \[\begin{array}{cc}\cos t&0\\0&1\end{array}\]}$.
So in this case $\det(F_1(t))=1>\cos t=\det(F_2(t))$ for $0<t<2\pi$.
For the second example again use $A_1=A_2=0$, but this time let
$R_1(t)\equiv {\scriptsize \[\begin{array}{cc}0&0\\0&1\end{array}\]}$
and  $R_2(t)\equiv 0$.  Then $R_1(t)\le R_2(t)$,
$F_1 ={\scriptsize \[\begin{array}{cc}1&0\\0&\cos t\end{array}\]}$,
and $F_2(t)={\scriptsize \[\begin{array}{cc}1&0\\0&1\end{array}\]}$.
So this time $\det(F_1(t))=\cos t<1=\det(F_2(t))$ and the inequality goes
in the other direction.

If $\la\,,\ra$ is positive definite and $A$, $B$ are self-adjoint with 
$A\le B$, then $\trace(A)\le\trace(B)$.  However if $\la\,,\ra$ is
indefinite then $A\le B$ does not imply any inequality between 
$\trace(A)$ and $\trace(B)$.  This is exactly where the proof of the volume
comparison theorem breaks down in the indefinite metric case.~\done
\end{remark}

\subsection{Comparisons for Tubes}

Theorem~\ref{hypercomp} implies a comparison theorem for the second
fundamental form of a  hypersurface parallel to a given hypersurface.
In the Riemannain case there are comparison results for tubes about
submanifolds of higher codimension which have interesting applications
(cf.~\cite{Heintze-Karcher}).  In the semi-Riemannain case there is
also a comparison result for the second fundamential forms of tubes
about non-degenerate submanifolds of $(\am,\ame)$, but there do not
seem to be applications of this result that have the same interest as
the Riemannian theorem.  Thus for the sake of completeness we include
the statement of the result, but omit the proof.  As with the
results above we state this as a theorem about systems of ordinary
differential equations.

Let $T$ be a subspace of $E$ and assume that the restriction of the
inner product $\la\,,\ra$ to $T$ is nondegenerate.  Let $T\oc$ be the
orthogonal complement of $T$ in $E$.  Then $E=T\oplus T\oc$.  Let 
$P:E\to T$ and $P\oc:E\to T\oc$ be the orthogonal projections.  This
includes the case when $T=\{0\}$, so  $P\oc=I$.

\begin{thm}\sl\label{gencomp}
For $i=1,2$, let $A_i:T\to T$ be a self-adjoint linear maps and $R_i(t)$ as
in Theorem~\ref{hypercomp}.  Then define $F_i(t)$ by
$$
F_i''(t)+R_i(t)F_i(t)=0,\quad F_i(0)=P,\quad F_i'(0)=-A_i P+P\oc
\,; \qquad i=1,2. 
$$
At the points where $F_i(t)$ is non-singular define 
$S_i(t)=-F_i'(t)F_i(t)^{-1}$.  As before the maps
$S_i(t)$ are self-adjoint at the points where it is defined.  
If $S_1(t)$ and $S_2(t)$ are defined on 
$(0,b]$ and 
$$
A_1\le A_2,\qquad R_1(t)\le R_2(t)\quad\mbox{for all $t\in (0,b]$},
$$
then $S_1(t)\le S_2(t)$.  If $S_1(b)=S_2(b)$ then $A_1=A_2$ and 
$S_1\equiv S_2$, $R_1\equiv R_2$ on $(0,b]$.~\done 
\end{thm}

\begin{remark}\normalshape
There is one case where this is of interest, and
that is when $M$ is a point $p$ in a Lorentzian manifold.  Let
$S_{-}(p):=\{u\in T(\am)_p:\ame(u,u)=-1\}$.  Then $S_-(p)$ is a
Riemannian manifold isometric to two disjoint copies of the hyperbolic 
space $\r^{n-1}_0(-1)$.  Define a map $f_r:S_-(p)\to\am$  by
$f_r(u)=\exp_p(ru)$.   If $(\am,\ame)$ is timelike geodesically
complete and the curvature satisfies $\la \ar(X,Y)Y,X\ra\le 0$ for all
pairs of vectors $X,Y$ spanning a timelike two plane, then for any
$X$ tangent to $S_-(p)$ at $u$
$$
\la f_{r*u}X,f_{r*u}X\ra\ge\frac1{r^2}\la X,X\ra.
$$
This result is due to Flaherty~\cite{Flaherty:Lorentz} where he used it to
prove a version of the Cartan-Hadamard theorem for Lorentzian
manifolds.  We note that this also follows from the result above as
the derivative of $f_r$ is given by 
$$
\la f_{r*u}X,f_{r*u}Y\ra=\la F_u(r)X,F_u(r)Y\ra
$$
and 
\begin{align*}
\frac{d}{dt}\la F_u(t)X,F_u(t)Y\ra&=2\la F_u'(t)X,F_u(t)X\ra\\
&=2\la F_u'(t)F_u(t)F_u(t)^{-1}X,F_u(t)X\ra\\
&=-2\la S_u(t)F_u(t)X,F_u(t)X\ra.
\end{align*}
Also for small $t$ it is not hard to see $F(t)=(1/t)I+O(1)$.
Therefore the comparison result for tubes can be used to prove the
inequality.~\done 
\end{remark}

\section{Rigidity in Flat Spaces}\label{sec:rigid-flat}

Let $\r^n_k=\r^n_k(0)$ be the flat simply connected space form of
index~$k$.  This is $\r^n$ with the semi-Riemannian metric
$\la\,,\ra=-\sum_{i=1}^kdx_i{}^2+\sum_{i=k+1}^ndx_i{}^2$.  Let $\eta$ be a
unit vector in $\r_k^n$, that is $\la\eta,\eta\ra=+1$ or $\la \eta,\eta\ra=-1$.
Let $\e=\la\eta,\eta\ra$, and let $\Lambda_\eta(\cdot)=\e\la\cdot,\eta\ra$.
Then $\Lambda_\eta$ is the linear functional on $\r^n_k$ that has the
orthogonal compliment of $\eta$ as kernel, and so that
$\Lambda_\eta(\eta)=1$. 
\begin{definition}\sl
A subset $B\subset \r^n_k$ has the
{\em compact intersection property} with respect to $\Lambda_\eta$ if 
and only if for each compact interval $[a,b]\subset \r$ the set 
$B\cap \{x:\Lambda_\eta(x)\in [a,b]\}$ has compact closure.~\done
\end{definition}

\begin{thm}\sl \label{flatrigid}
Let $(\am,\ame)$ be a semi-Riemannian manifold of
dimension~$n\ge3$ so that: 
\begin{enumerate}
\item The curvature tensor of $\ame$ satisfies 
$\e\ar \geq 0$ on $\am$ in the sense of Definition~\ref{def:curvbound}.
\item Every geodesic of sign $\e$ is complete.  (That is if
$\gamma:(a,b)\to\am$ is a geodesic with 
$\ame(\gamma',\gamma')=\e$, then $\gamma$ extends to a geodesic defined on
all of $\r$.) 
\end{enumerate}
Let $\eta_0$ be a unit vector in $\r_k^n$ of sign $\e$, and 
$B\subset\r^n_k$ be a closed subset so that $\r_k^n\setminus B$ is
connected and 
\begin{enumerate}\setcounter{enumi}{2}
\item $B\subset \{x:\Lambda_{\eta_0}(x)>0\}$ 
\item $B$ has the compact intersection property with respect to 
$\Lambda_{\eta_0}$.
\end{enumerate}
Then any local isometry $\phi:\r^n_k\setminus B\to \am$ extents to a
surjective local
isometry $\widehat{\phi}:\r_k^n\to \am$ defined on all of $\r_k^n$.
\end{thm}

\begin{remark}\normalshape
For  Riemannian manifolds this does not lead to
any results in the case of sectional curvatures $\le0$.  However we
note that a proof by Schroeder and
Ziller~\cite[Theorem~1]{Schroeder-Ziller:Rigidity} of a result of
Gromov~\cite[\S 5]{ballmann:gromov:schroeder} leads to
\smallskip

\noindent{\bf Theorem\ }{\sl Let $(\am,\ame)$ be a complete simply
connected Riemannian manifold of dimension at least $n\ge2$ and 
non-positive sectional curvature.  Let $m\ge n$ and 
$B\subset \r^m$ have the compact intersection property with respect to
$\Lambda_\eta$ and assume $\r^m\setminus B$ is connected.  Then any
isometric imbedding $\phi:\r^n\setminus B\to \am$ as a totally
geodesic submanifold extends uniquely to
a isometric imbedding $\widehat{\phi}:\r^m\to\am$ and a totally
geodesic submanifold.}~\done
\smallskip

\noindent This is interesting in that it is not assumed
$B\subset \{x:\Lambda_\eta(x)>0\}$ or that the model space $\r^m$
and the space $\am$ have the same dimension.  On the other hand it is
important in their proof (which is based on Toponogov's triangle
comparison theorem) that the map $\phi$ be injective.

When $(\am,\ame)$ is Riemannian, the sectional curvatures of
$(\am,\ame)$ are $\ge 0$, the set $B\subset \r^n$ is assumed
compact and the map $\phi:\r^n\setminus B \to  \am$ 
is injective, the result
can be deduced from a result of Greene and  
Wu~\cite[Theorem~1]{Greene-Wu:gap-theorems} (see also
\cite[Remark on p.~75]{ballmann:gromov:schroeder}). When $\phi$ is injective
it is easy to see that the growth rate of the volume of geodesic balls
is the same as that of balls in Euclidean space and so the
Bishop-Gromov volume comparison theorem can be used to prove the
result under the weaker assumption that the Ricci tensor of
$(\am,\ame)$ is non-negative, cf.~\cite{Karcher:comp}.
If $(\am,\ame)$ is a spin manifold and $\am\setminus\phi[\r^n\setminus
C]$ is compact then by rigidity results related to the positive mass
conjecture it is enough to assume that the scalar curvature of
$(\am,\ame)$ is non-negative (cf.~\cite{bartnik:mass}).
If the dimension of $\am$ is $\le 6$
then the proofs of Schoen and Yau
\cite{schoen:yau:ADM,schoen-yau:scalar} of the positive mass
conjecture
imply the result without the
assumption that $(\am,\ame)$ is spin.
For other related rigidity results in the Riemannian case
cf.~\cite{Kasue-Sugahara:Gap} and \cite{Greene-Petersen-Zhu:decay}.~\done
\end{remark}

The next corollary is not much more than a special case of Theorem 
\ref{flatrigid}, but in doing the inductive step of the proof of the theorem 
and as a lemma for use in later sections, it is helpful
to have it stated separately.

\begin{cor}\sl \label{domain}
Let $D\subset \r_k^n$ be a connected domain with compact
closure and smooth boundary and $\ame$ be a semi-Riemannian metric 
on $D$ so that $\ame$ agrees with the standard metric $g_0$ in a 
neighborhood of $\partial D$. 
If the curvature tensor satisfies either
$\ar \ge 0$ on all of $D$, or $\ar \le 0$ 
on all of $D$.  Then $\ar\equiv 0$.
\end{cor}

\begin{pf}  Define a new semi-Riemannian metric $\ame_1$ on $\r^n$ by 
letting $\ame_1=\ame$ in $D$ and $\ame_1=g_0$ in $\r^n\setminus D$.
As $\ame$ and $g_0$ agree in a neighborhood of $\partial D$ this metric is
smooth.  We first consider the case  $\ame$, and thus $\ame_1$, is
positive definite.  The hypotheses imply  the sectional curvature
of $\ame_1$ is either non-negative or non-positive.  In these cases
theorems of Greene and Wu~\cite[Theorem 1]{Greene-Wu:gap-theorems}
(for the case $\ar\ge0$) and  Kasue and 
Sugahara~\cite[Theorem~2]{Kasue-Sugahara:Gap} (for the case $\ar\le 0$)
imply that $\ar \equiv0$.  
The case that $\ame$ is negative definite reduces 
to the  positive definite
case by replacing $\ame$ by $-\ame$.  In all other cases it will be
possible to choose a unit vector $\eta$ so that
$\la\eta,\eta\ra\ar \geq 0$.  Then since $D$ is compact it will
have the compact intersection property with respect to $\Lambda_\eta$.
So this result follows from Theorem~\ref{flatrigid}.
\end{pf}

\begin{pf*}{Proof of Theorem~\ref{flatrigid}} The proof is by induction on the 
dimension $n$ of $\am$.  The base of the induction is
Theorem~\ref{thm:twodrigid}.  Now some notation is needed.  Let
$N(\eta_0)=\{x\in \r_k^n:\Lambda_{\eta_0}(x)=0\}$ be the hyperplane through
the origin, 
orthogonal to $\eta_0$, and set 
$M(\eta_0):=\phi[N(\eta_0)]$.  
Let $\nor_0$ be the unit normal field to $M(\eta_0)$
so that $\nor_0(\phi(0))=\phi_{*0}\eta_0$.  For $r\in \r$ define
$f_{\eta_0,r}:N(\eta_0)\to \am$ by 
$$
f_{\eta_0,r}(x)=\exp_{\phi(x)}(r\nor_0(\phi(x))) ,
$$
where $\exp$ is the exponential map of the metric $\ame$.  Thus
$M(\eta_0)[r]:=f_{\eta_0,r}[M(\eta_0)]$ is the parallel hypersurface 
to $M(\eta_0)$ at a distance $r$.

Our next goal is to show  the curvature tensor $\ar$ of $(\am,\ame)$
is zero at any point of the form $f_{\eta_{0},r}(x)$.  If $x\in N(\eta_0)$
and the segment $t\mapsto x+t\eta_0\notin B$ ($0\le t\le r$) then, 
as $\phi$ is a local isometry on $\r^n_k\setminus B$, we have 
$\ar_{f_{\eta_0,r}(x)}=0$, and also the image of $f_{\eta_0,r}$ near 
$f_{\eta_0,r}(x)$ is given by $\phi \{y\in \r_k^n:\Lambda_{\eta_0}(y)=r\}$.
Therefore near $f_{\eta_0,r}(x)$ the image of $f_{\eta_0,r}$ is totally
geodesic and $f_{\eta,r}:N(\eta_0)\to M(\eta_0)[r]$ is a local isometry
near $x$.  Set
$$
r_0=\sup\{r: \ar_{f_{\eta_0,s}(x)}=0\ \mbox{\rm{for all}}\ x\in N(\eta_0)\ 
\mbox{\rm and}\ s\le r\}.
$$

If $r_0=\infty$, then $\ar_{f_{\eta_0,r}(x)}=0$ for all $x\in N(\eta_0)$
and $r\in \r$ as claimed.  Thus assume, toward a contradiction,  
$r_0<\infty$.  From the hypothesis of the theorem $r_0>0$.  If $r\le r_0$
the above discussion shows  $f_{\eta_0,r}:N(\eta_0)\to M(\eta_0)[r]$
is a local isometry and so it is an immersion.  By continuity if 
$r>r_0$ is only slightly larger than $r_0$ then 
$f_{\eta_0,r}:N(r_0)\to \am$ will also be an immersion.  For such an $r$
consider the hypersurface $M(\eta_0)[r]=f_{\eta_0,r}[N(\eta_0)]$.  
By the compact intersection property there is a connected domain with 
smooth boundary $D_0\subset N(\eta_0)$ with compact closure and so that 
if $x\in N(\eta_0)\setminus D_0$ and $s\le r$ then $x+s\eta_0\notin B$.
Therefore the pulled back metric $f^*_{\eta_0,r} \ame$ is the usual flat
metric on $N(\eta_0)\setminus D_0$.  

We now use the comparison theorem.
Using the Riccati equation satisfied by the Weingarten map of a family of
parallel hypersurfaces (Proposition~\ref{2ndForm}) and one of the comparison 
results for Riccati equations (Theorem~\ref{curvcomp}), if
$S_{f_{\eta_0,r}(x)}$  is the Weingarten map of $M(\eta)[r]$ at 
$f_{\eta_0,r}(x)$ then $\W^2(S_{f_{\eta_0,r}(x)})\ge 0$.  Let 
$R^{\eta_0,r}$ be the curvature tensor of $M(\eta_0)[r]$.  Then by the
Gauss curvature equation (\ref{gauss}) multiplied by $\e$ and the
assumptions on $\ar$, we have
$$
\e \la R^{\eta_0,r}(X,Y)Y,X\ra
=\e \la \ar(X,Y)Y,X\ra+ \la {\W}^2(S_{f_{\eta_0,r}(x)})X\w Y,X\w Y\ra\ge 0
$$
for any $X,Y$ tangent to $M(\eta)[r]$.
Therefore the metric $f_{\eta_0,r}{}^*\ame$ agrees with the usual metric on
$N(\eta_0)$ outside of $D_0$ and the curvature satisfies 
$\e\la R^{\eta_0,r}(X,Y)Y,X\ra\ge 0$.  
Thus by the induction hypothesis this
implies $R^{\eta_0,r}\equiv 0$.  For this to hold we must have 
$\la \W^2(S_{f_{\eta_0,r}(x)})X\w Y,X\w Y\ra=0$ for all $X$, $Y$ tangent to 
$M(\eta)[r]$.  Putting this in the last equation implies 
$\la \ar(X,Y)Y,X\ra=0$ for all $X$, $Y$ tangent to $M(\eta_0)[r]$.
Unfortunately this is not enough to conclude directly that full curvature 
tensor $\ar$ vanishes along $M(\eta_0)[r]$.  

Let $y_1=f_{\eta_0,r_1}(x_1)$ be so that $\ar_{y_1}\ne 0$.  
Assume  
$r_1$ is taken to only be slightly larger than $r_0$ so that $f_{\eta_0,r_1}$
is an immersion.  Let $\eta\in\r_k^n$ be a unit vector.  Then define
$N(\eta)=\{ x\in \r_k^n:\la x,\eta\ra=0\}$, etc.\ 
just as was done in the case 
of $N(\eta_0)$.  Note however that in general $N(\eta)\cap
B\ne\emptyset$ and so $M(\eta)=\phi[N(\eta)\setminus B]$ and the maps
$f_{\eta,r}$ are only defined on $N(\eta)\setminus B$. 
Let $D_1\subset N(\eta_0)$ be a connected domain with smooth boundary
so that for all $r\le r_1$, $B\cap (N(\eta_0)+r\eta_0)$ is a subset of the
interiors of $D_1+r\eta_1$.  Then the pulled back metric
$f^*_{\eta_0,r_1}\ame$ agrees with the standard flat metric of $N(\eta_0)$
on $\partial D_1$.  Now let $\eta$ be a unit vector of $\r_k^n$ which is
very close to $\eta_0$. Then for some $r$ very close to $r_1$ there holds
$y_1=f_{\eta,r}(x)$ for some $x\in N(\eta)$. Note here that $r,x$ are
functions of $\eta$.

If $\eta$ is close 
enough to $\eta_0$ there is a 
domain $D\subset N(\eta)$, which can be taken to be close to the domain
$D_1\subset N(\eta_0)$, so that the pulled back metric $f^*_{\eta,r}\ame$
agrees with the standard flat metric of $N(\eta)$ on a neighborhood of 
$\partial D$.  The comparison argument of the last paragraph still holds
and so the curvature tensor $R^{\eta,r}$ of $f^*_{\eta,r}\ame$ satisfies
$\e\la R^{\eta,r}(X,Y)Y,X\ra\ge 0$.  Thus by the induction hypothesis and
Corollary~\ref{domain}, $R^{\eta,r}\equiv 0$ on $D$.  Again this implies
$\la \ar(X,Y)Y,X\ra=0$ for all vectors $X$, $Y$ tangent to 
$f_{\eta,r}[D]$.

Let $\xi(\eta)$ be the normal to $f_{\eta,r(\eta)}[D]$ at $y_1$.  Then the map
$\eta\mapsto \xi(\eta) \in T(\am)_{y_1}$ is smooth in a neighborhood 
of $\eta_0$, and if 
we choose $r_1$ close enough to $r_0$ we also have that the derivative
$\xi_*$ is non-singular at $\eta_0$.  By the implicit function theorem this
implies  there is a neighborhood $U$ of $\xi(\eta_0)$ in the tangent space
$T(\am)_{y_1}$ so that if $u$ is a unit vector in $U$, then there is a 
unit vector $\eta$ near $\eta_0$ with $\xi(\eta)=u$.  This in turn implies
that if $X_0$, $Y_0$ are linearly independent vectors in
$T(f_{\eta_0,r_1}[N(\eta_0)])_{y_1}=\xi(\eta_0)\oc$ 
then there are neighborhoods
$V$ of $X_0$ and $W$ of $Y_0$ in $T(\am)_{f_{y_1}}$ so that if $X\in V$ and 
$Y\in W$ then there is a $\eta$ near $\eta_0$ so that both $X$ and $Y$ are 
orthogonal to $\xi(\eta)$.  Thus $X,Y\in \xi(\eta)\oc=T(f_{\eta,r}[D])$.
Therefore by the discussion of the last paragraph $\la\ar(X,Y)Y,X\ra$=0.
The map $(X,Y)\mapsto \la\ar(X,Y)Y,X\ra$ is a polynomial map on
$T(\am)_{y_1}\times T(\am)_{y_1}$, and this map vanishes on the open set
$V\times W$.  Thus $\ar_{y_1}=0$.  However the point $y_1$ was chosen so
that $\ar_{y_{1}}\ne 0$.  This contradiction completes the proof of the
earlier claim  that $\ar=0$ at any point of the form $y=f_{\eta_0,r}$.

Note that points of $\r_k^n$ have a unique expression of the form
$x+r\eta_0$ where $x\in N(\eta_0)=\eta_0\oc$ and $r\in\r$.  Define  
$\widehat{\phi}:\r^n_k\to \am$ by
$$
\widehat{\phi}(x+r\eta_0)=\exp_{\phi(x)}(r\nor_0(\phi(x)))=f_{\eta_0,r}(x).
$$
Using that $\ar=0$ at all points of the form $f_{\eta_0,r}(x)$ it is not
hard to show  $\widehat{\phi}$ is a local isometry.
This gives the required extension of $\phi$ to $\r_k^n$.  That
$\widehat{\phi}$ is surjective follows from Proposition~\ref{onto}.
\end{pf*}

\section{Rigidity in Simply Connected Space Forms}\label{sec:space-forms}

The results in this section are for Riemannian and Lorentzian
manifolds. The central ideas are an inductive argument using parallel
hypersurfaces and rigidity results for warped products 
that follow directly from the comparison theory. 

\subsection{Rigidity Lemmas in Warped Products}\label{sec:warpsplit}

Let $(M,g)$ be a Riemannian manifold and let $w:[0,L]\to
(0,\infty)$ be a positive $C^2$ function with $w(0)=1$.  Set
$\rr =-w''/w$ and $a=w'(0)$.  Then
$$
w''(t)+\rr(t)w(t)=0, \quad w(0)=1, \quad w'(0)=a.
$$
Let $\e=+1$ or $\e=-1$ and let $\bargw$ be the warped product metric
$$
\bargw:=w(t)^2g+\e dt^2
$$
on $\am:=M\times[0,L]$.  Then for each $x\in M$ the curve
$c_x(t):=(x,t)$ is a unit speed geodesic in $(M\times[0,L],\bargw)$.
Then $\T$ is a normal to the hypersurface $M\times\{t\}$ and with
respect to this normal the Weingarten map is $S(t)=-w'(t)/w(t)I$.  In
particular the Weingarten map of $M\times\{0\}$ is $-w'(0)I=aI$.
\begin{definition} \sl\label{def:prod-apt}
A semi-Riemannian metric $\ame$ on $M\times[0,L]$ is {\sl adapted to the
product structure\/} iff for each $x\in M$ the curve $c_x(t)$ is a
unit speed geodesic with respect to the metric $\ame$.~\done  
\end{definition}
If $x^1\cd x^{n-1}$ are local coordinates on $M$, then $x^1\cd
x^{n-1},t$ are local coordinates on $M\times [0,L]$ and the metric
$\ame$ on $M\times[0,L]$ is adapted to the product structure if and
only if the metric is of the form 
$\ame=\sum_{i,j=1}^{n-1}g_{ij}(x^1\cd x^{n-1},t)dx^idx^j+\e dt^2$.  
\begin{remark}\normalshape
\label{rmk:gauss-adpt}
Let $f:(M,g)\to(\am,\ame)$ be an isometric immersion of 
$(M,g)$ into $(\am,\ame)$ as a hypersurface there is a globally
defined unit normal $\eta$ along $f$.  Define a map $F:M\times[0,L]\to
\am$  by $F(x,t)=\exp(t\eta(x))$.  If this is a local diffeomorphism
then the Gauss lemma implies the metric $F^*\ame$ is adapted to the
product structure of $M\times[0,L]$.~\done
\end{remark}

If $\ame$ is adapted to the product structure of $M\times[0,L]$ then
$\T$ is a unit normal to the hypersurfaces $M\times\{t\}$.  
Denote by
$\aS(t):=\ac_X\T$ the Weingarten map of $M\times\{t\}$ in
the metric $\ame$.  Along each of the geodesics $c_x(t)=(t)$ there 
and let $\ar(t)=\ar(\cdot,\T)\T$.  

\begin{thm}\sl\label{warprigid1} 
With the notation above (so
that the restriction of the metric $\ame$ to $M\times\{0\}$ is positive
definite),  assume that one of the two conditions 
\begin{enumerate}
\item $S(0)\ge aI$ on $M\times\{0\}$ and $\ar(t)\ge \rr(t)I$ 
on $M\times [0,L]$ or,
\item $S(0)\le aI$ on $M\times\{0\}$ and $\ar(t)\le \rr(t)I$ 
on $M\times [0,L]$.
\end{enumerate}
holds.  If
\begin{equation}\label{equiv}
\aS(L) =-\frac{w'(L)}{w(L)}I\quad\mbox{on all of}\quad M\times\{L\}
\end{equation}
then $\ame\equiv\bargw$.
\end{thm}

\begin{pf} Assume the first condition holds. The equality (\ref{equiv}) implies
that equality holds between $S_1(t):=-\frac{w'(t)}{w(t)}I$ and 
$S_2(t):=\aS(t)$ at $t=L$ in the comparison Theorem~\ref{hypercomp}. 
This implies that 
$S_1\equiv S_2$ and $\rr(t)I\equiv \ar(t)$.  The rest follows  by a
direct calculation.
\end{pf}

\begin{remark}\normalshape For Theorem \ref{hypercomp} to apply we 
need that the identity map on tangent spaces to $(M,g)$ is positive definite. 
This will only be the case when $g$ is positive definite, so that 
$(\am, \ame)$ will be Riemannian or Lorentzian. This is the reason why the 
results of this section are restricted to these cases.~\done
\end{remark}
\begin{thm}\sl\label{warprigid2}
With the notation above assume that one of the conditions
\begin{enumerate}
\item $S(0)= aI$ on $M\times\{0\}$ and $\ar(t)\ge \rr(t)I$
 on $M\times [0,L]$,and $-\frac{w'(r_0)}{w(r_0)}>0$,
\item $S(0)= aI$ on $M\times\{0\}$ and $\ar(t)\le \rr(t)I$ 
on $M\times [0,L]$,
and $-\frac{w'(r_0)}{w(L)}<0$,
\item $S(0)= aI$ on $M=M\times\{0\}$, $\ar(t)\le \rr(t)I$ 
and $\aS(t)>0$ on $M\times [0,L]$, or
\item $S(0)= aI$ on $M\times\{0\}$, $\ar(t)\ge \rr(t)I$ 
and $\aS(t)<0$ on 
$M\times [0,L]$
\end{enumerate}
holds.  Then 
$$
\W{}^2(S(L))= \(\frac{w'(L)}{w(L)}\)^2\W{}^2(I)\quad\mbox{on all of}
\quad M\times\{L\}
$$
implies $\ame\equiv\bargw$.
\end{thm}

\begin{pf} This follows from the comparison Theorems~\ref{defcompcurv} and 
\ref{defcompcurv2} in the
same way that the last theorem follows from Theorem~\ref{hypercomp}.
\end{pf}

\begin{remark}\normalshape
\label{curvremark}
If the curvature tensor of $(\am,\ame)$ satisfies
$$
\ame(\ar(X,Y)Y,X)\ge \pm \e\rr(t)(\ame(X,X)\ame(Y,Y)-\ame(X,Y)^2)
$$
on all of the image of $G:M\times[0,r_0]\to \am$ then $\ar(t)\ge \pm \rr(t)I$.
Likewise
$$
\ame(\ar(X,Y)Y,X)\le \pm\e\rr(t)(\ame(X,X)\ame(Y,Y)-\ame(X,Y)^2)
$$ 
on the image of $G:M\times[0,r_0]\to \am$ implies 
$\ar(t)\le \pm\rr(t)I$.~\done
\end{remark}

\subsection{The Rigidity Results}\label{sec:rigid}

\begin{definition}\sl Let $N$ be a manifold and $\tau:N\to\r$ be a
function.  Then $B\subset N$ has the {\em compact intersection
property\/} with respect to $\tau$ if for every compact interval
$[a,b]\subset\r$ the set $\{x\in B:\tau(x)\in [a,b]\}$ has compact closure
in $N$.~\done
\end{definition}

The rigidity results of this section are of the following type:  Let
$k =0$ or $k = 1$ and let 
$(\mod,\mme)$ be a model space (by which we mean one of the spaces
listed in Table~2) and $(\am,\ame)$ a semi-Riemannian manifold of
the same dimension and index as $\mod$.  Let $B\subset\mod$ be a
closed set and $f:(\mod\setminus B)\to  \am$ a local isometry.  Then
we  are looking for conditions on $B$ and the curvature tensor of
$(\am,\ame)$ so that $f$ will extend to a local isometry defined on
all of $\mod$.  The basic conditions are that $B$ not be too large
(which will usually mean that it have the compact intersection property
with respect to some function $\tau$), that $\mod\setminus B$ be connected
(which is easily seen to be necessary) and that $(\am,\ame)$ is
geodesically complete with respect to geodesics of some sign. 
Table \ref{table2} summarizes the conditions on the model space, the curvature 
bounds, the sign of the complete geodesics (column headed by Geo.), the 
exceptional set, and the function.

\begin{table}[htp]
\caption[]{}
\label{table2}
\begin{center}
\begin{tabular}{|c|c|c|c|c|}
\hline
Model & Curv. & Geo.&Exceptional Set & Function  \\
\hline\hline
$S^n=\r_0^n(+1)$ & $\stack{\vs\ar\geq +1}{\ar\leq +1}$ &$+$& $B \subset S^n\setminus S^{n-1}$    &    ---     \\
\hline
$\r^n_0(-1)$       & $\stack{\vs\ar\geq -1}{\ar\leq -1}$ &+& $B \subset \{ \tau > 0\}$   &   $\tau = w=e^{-t}$    \\
\hline
$\r^{n}_1(+1)$ & $\stack{\vs\ar\geq +1}{\ar\leq +1}$&$-$& $B\subset \{ \tau >0\}\subset\h_1^n(+1)$& $\tau = w=e^{t}$ \\
\hline
$\r^{n}_1(+1)$ & $\stack{\vs\ar\geq +1}{\ar\leq +1}$&$-$& $B \subset \r^n_1(+1) \setminus S^{n-1}$&     ---      \\
\hline
$\r^{n}_1(-1)$ & $\stack{\ar\geq -1}{\ar\leq -1}$&$-$& $B \subset \{ \tau > 0\}$ &\begin{minipage}{.75in}$\tau=$ time\\
						function of\\ a geodesic\end{minipage}    \\
\hline
$\r^{n}_1(-1)$ & $\stack{\vs\ar\geq -1}{\ar\leq -1}$&$-$& $B\subset \{\tau>0\}\subset{\bf S}^{n}_1(-1)$&$\tau=w=\cosh(t)$\\
\hline
\end{tabular}
\end{center}
\end{table}

In the first two rows of the table $S^{n-1}$ is imbedded in
$S^n=\r^n_0(+1)$ as an equator (that is as a totally geodesic
submanifold).  In the next two rows $\r^n_0(-1)$ is viewed as a warped
product as in row~1 of Table~1 of Section~\ref{sec:warpmodel} and the
function $\tau$ is the warping function $w$.  Likewise in the next two
rows $B$ is a subset of the ``half space'' $\h^n_1(+1)$ in the 
de~Sitter space $\r^n_1(+1)$ which is represented as a warped product
in row~2 of Table~1 and in this case the function $\tau$ is again
the warping function (see also Definition~\ref{def:half-space}).  
The de~Sitter space $\r^n_1(+1)$ has a
another representation  as a warped product $-dt^2+\cosh^2(t)\ame_{{\bold
R}^{n-1}_0(+1)}$ and thus $S^{n-1}=\r^{n-1}_0(+1)$ is being viewed as
the submanifold $S^{n-1}\times\{0\}$ of $\r^n_1(-1)$.  In the first
two rows for the anti-de~Sitter space $\r^n_1(-1)$, the function
$\tau$ is the time function of a timelike unit speed geodesic as in
Definition~\ref{def:time-fcn}.  Finally in the last two rows of the table
$B$ is a subset of ``strip'' $\els_1^n(-1)$ in the 
anti--de~Sitter space as in Definition~\ref{def:strip} and the function
$\tau$ is the warping function. 

\begin{thm}\sl \label{thm:rigidity-curved}
Let $(\mod,\mme)$ be one of the model spaces in
Table~2 of dimension at least three
and let $B\subset\mod$ be a closed set so that $\mod\setminus
B$ is connected.  Then for any of the twelve rows in the table if $B$
is a subset of the indicated set and if 
a function $\tau$ is given in the last column assume that $B$
has the compact intersection property with respect to $\tau$.  Let
$(\am,\ame)$ be a semi-Riemannian manifold of the same dimension and
index as $(\mod,\mme)$ so that the curvature tensor $\ar$ of
$(\am,\ame)$ satisfies the indicated inequality and that every geodesic
of $(\am,\ame)$ of the indicated sign is complete.  Then any local
isometry $\phi:\mod\setminus B\to \am$ extends to a surjective local
isometry $\whp:\mod\to \am$.  Thus $(\am,\ame)$ has constant
curvature. 
\end{thm}

\begin{remark}\normalshape 
(1)
In the case of $\r^n_0(-1)$ (the Riemannian hyperbolic space), when
$B$ is compact,  $\phi:\r^n_1(-1)\setminus B\to \am$ is injective, and
$\am$ is simply connected
(and when $\ar\ge0$ also assume that $(\am,\ame)$ has a pole, that is a
point where $\exp_{x_0}:T(\am)_{x_0}\to \am$ is a diffeomorphism) the
result can be deduced from a result of Kasue and 
Sugahara~\cite[page~697]{Kasue-Sugahara:Gap}.  If $\ar\le -1$, $\am$
is simply connected and $\am\setminus\phi[\r_0^n(-1)\setminus B]$ is
compact this is also a special case of a result of Schroeder and
Ziller~\cite[Theorem~7]{Schroeder-Ziller:Rigidity}.  We note that our
result allows the set $B$ to be non-compact and does not require
$\phi$ to be injective.

For the space $\r^n_0(+1)=S^n$ when $\ar\ge+1$ if the space
$\am\setminus\phi[S^n\setminus B]$ is simply connected, compact and
strictly convex then the result is covered by the result of Schroeder and
Ziller~\cite[Theorem~7]{Schroeder-Ziller:Rigidity}. In the case of
$\r^n_0(+1)=S^n$ and $\ar\le +1$ to the best of our knowledge the
result is new.

(2) Because of the failure of Theorems \ref{warprigid1} and \ref{warprigid2}
when $(\am, \ame)$ is not Riemannian or Lorentzian, the method of proof here 
fails as a method for proving rigidity theorems for semi-Riemannian manifolds 
of arbitrary index.~\done
\end{remark} 

The proofs all follow the same pattern: An induction on dimension 
using hypersurfaces parallel to a standard hypersurface and the warped
product rigidity results of Section~\ref{sec:warpsplit}.
We will do one case in detail and leave the
others to the reader.  To avoid inessential but annoying problems
involving signs, the case we consider is the Riemannian case where all the
geometric complications involved occur.  

\begin{thm}\sl\label{thm:NewSphereThm}
Let $(\am,\ame)$ be a complete Riemannian manifold of dimension
$n\ge3$ with sectional
curvatures $\le1$.  Let $B\subset S^n\setminus B$ be a closed set with
$S^n\setminus B$ connected and let $\phi :S^n\setminus B \to \am$
be a local isometry.  Then $\phi$ extends to a surjective local
isometry $\widehat{\phi}:S^n\to \am$.
\end{thm}

In doing the induction step it is useful to have the following special
case.  The notation is as follows:
$(S^{n-1}\times(-r_0,r_0),\elm)$ is the tube of radius $r_0$ about the
equator $S^{n-1}$.  
\begin{cor}\sl\label{sphererigid2}
Let $(\am,\ame)$ be a connected complete Riemannian manifold of dimension 
at least three and so that the sectional curvature of $(\am,\ame)$ satisfies 
$K_{\am}\le 1$.  Let $0<r_0<\pi/2$.  If $(\am,\ame)$ contains a subset
isometric to $(S^{n-1}\times(-r_0,r_0),\elm)$, 
then $(\am,\ame)$ is isometric to $(S^{n},\elm)$.
\end{cor}

\begin{pf} By assumption there is an isometry
$\phi:(S^{n-1}\times(-r_0,r_0),\elm)\to(\am,\ame)$.  By the theorem this
extends to a local isometry $\whp:(S^{n},\elm) \to  (\am,\ame)$.  
As $\whp$ is
a local isometry the image $\whp[S^n]$ is open in $\am$.  But as $S^n$
is compact the image is also closed.  As $\am $ is connected this
implies  $\whp$ is surjective.  To show   $\whp$ is an isometry it
only remains to show  it is injective.

Note  $\whp:S^{n}\to\am$ is a
covering map.  Let $G$ be the group of deck transformations of 
$\whp:S^{n}\to\am$, that is $G$ is the set of maps $a:S^n\to S^n$ so
that $\whp\circ a=\whp$.  If $a\in G$, then $a:S^n\to S^n$ is an
isometry of $\elm$.  It follows  
$a[S^{n-1}]\cap S^{n-1}\ne \emptyset$.  Let
$x\in S^{n-1}$ so that $a(x)\in S^{n-1}$.  Then $\whp(ax)=\whp(x)$.  But 
on $S^{n-1}$ we have $\whp=\phi$ so $\phi(ax)=\phi(x)$. 
But $\phi$ is an isometry
and therefore injective.  Thus $a(x)=x$.  But a deck transformation with a
fixed point is the identity so $a=\mathop{\mathrm{Identity}}$.  As $a$ was an
arbitrary element of $G$ this implies $G$ is trivial.  As $S^n$ is simply
connected this implies $\whp$ is injective and completes the proof 
$(\am,\ame)$ is isometric to $(S^n,\elm)$.
\end{pf}

\begin{pf*}{Proof of Theorem~\ref{thm:NewSphereThm}} 
It is enough to show $(\am,\ame)$ has constant sectional curvature~1.
For then the universal covering space of $(\am,\ame)$ is the standard
sphere $(S^n,\elm)$.  Let $\pi:S^n\to\am$ be the covering map.  As any
local isometry between connected open sets in $S^n$ extends to a
global isometry there is an isometry $\psi:S^n\to S^n$ so that
$\pi\circ\psi|_{S^n\setminus B}=\phi$.  
Then $\widehat{\phi}:=\pi\circ\psi$ is the
required extension. That $\widehat{\phi}$ is surjective follows from 
Proposition~\ref{onto}.

The rest of the proof is devoted to
showing that $(\am,\ame)$ has constant sectional curvature~$1$.
Let $B$ be the set of points $x$ of $\am$ where some sectional curvature at 
$x$ is less than $+1$.  
Toward a contradiction assume 
$B\ne\emptyset$.  We may assume  $r_0$ is maximal with respect to the
property that there is a local isometry 
$\phi:(S^{n-1}\times[-r_0,r_0],\elm)\to (\am,\ame)$. For $U\subset \am$ 
let $\cl (U)$ be the closure of $U$ in $\am$.  Then as $r_0$ is maximal 
there is a point $y_0\in \phi[S^{n-1}\times[-r_0,r_0]]\cap \cl(B)$.  
By acting with an isometry on $S^{n-1}\times\{0\}$ we may choose the 
foliation defined by $S^{n-1}\times\{r\}$ so that in fact 
$y_0 = \phi[S^{n-1}\times[-r_0,r_0]]\cap \cl(B)$.  
Let 
$x_0\in S^{n-1}\times[-r_0,r_0]$ so that $\phi(x_0)=y_0$.  We may assume
 $x_0\in S^{n-1}\times\{r_0\}$. If $x_0\in S^{n-1}\times\{-r_0\}$ 
then do the change of variable $(x,t)\mapsto (x,-t)$ on
$S^{n-1}\times[-r_0,r_0]$.

Let $K_0=1/\cos^2(r_0)$ so that $K_0$ is the sectional curvature of 
$S^{n-1}\times\{r_0\}$ with the metric induced by $\elm$.  Let 
$0<r_1<r_0$ and let $K_1=1/\cos^2(r_1)$.

Let $S^{n-1}(K_1)$ be the
totally umbilic submanifold of $(S^{n-1}\times[-r_0,r_1])$  which is
tangent to $S^{n-1}\times\{r_0\}$ at the point $x_0$.  Then 
$S^{n-1}(K_1)\cap S^{n-1}\times\{r_0\}=\{x_0\}$ and 
the Weingarten map of $S^{n-1}(K_1)$ is 
$S\equiv (\sin(r_1)/\cos(r_1))I$ (cf.~Table~1).  If $r_1$ is chosen
close enough to 
$r_0$ then $S^{n-1}(K_1)\subset S^{n-1}\times[0,r_0]$.

Let $\eta$ be the unit normal along $S^{n-1}(K_1)$ so that
$\eta(x_0)=\T{}_{x_0}$.
Let $r_3$ be a small positive number to be chosen later and define 
a map $G:S^{n-1}(K_1)\times [0,r_3]\to \am$ by
$$
G(x,t)=\exp_{\phi(x)}(\phi_*\eta(x)).
$$
As $\phi$ is a local isometry on $S^{n-1}\times [-r_0,r_0]$ the submanifold
$\phi[S^{n-1}(K_1)]$ is also totally umbilic in $\am$ with Weingarten map
$(\sin(r_1)/\cos(r_1))I$.  
For $t>0$ let $S_t$ be the Weingarten map of 
$S^{n-1}(K_1) \times\{t\}$ in $S^{n-1}(K_1) \times[0,r_3]$ 
with respect to the 
pullback metric $G^*\am$.  
If $u(t):=\tan(r_1 +t)$ then $u'=u^2+1$.  
By the assumption on the curvature and the comparison Theorem
\ref{defcompcurv2},
$$
\W{}^2(S_t)\le u(t)^2\W{}^2(I)\quad \mbox{on}\ [0,r_3] 
$$
This, the assumption that the sectional curvature of $\ame$ (and thus also
the pulled back metric $G^*\ame$) has sectional curvature $\le 1$ implies
the sectional curvature of the hypersurface $S^{n-1}(K_1) \times\{r_3\}$ with
the metric induced by $G^*\ame$ satisfies 
\begin{equation}
K_{S^{n-1}(K_1)\times\{r_3\}}\le \frac{1}{\cos^2(r_1+r_3)}.  \label{crap}
\end{equation}

Let $K_3:=1/\cos^2(r_1+r_3)$. 
We now claim  if $r_3>0$ is small enough, then equality holds in 
inequality~(\ref{crap}).  This follows from the theorem we are proving by
induction.  The base case is $n=3$ so that $S^{n-1}(K_1)\times \{r_3\}$ 
is two
dimensional.  Let
$$
\rho(r_3):=1/\cos(r_1+r_3) .
$$
By making $r_3$ small enough, the Gauss curvature $K$ of the
rescaled metric $(S^{2}(K_1)\times \{r_3\} ,\rho(r_3)G^* \ame)$
satisfies $K\le1$ and can be made arbitrarily close to $1$.

Also by making $r_3$ small we can make the set where $K=1$ as 
large as we please, in the sense that there is a subset of 
$(S^{2}(K_1)\times \{r_3\},\rho(r_3)^2  G^*\ame)$ which is 
isometric to a standard sphere
of constant Gauss curvature~$+1$ with a very small ball deleted.  Thus the 
two dimensional rigidity result Theorem~\ref{thm:calabi} implies 
$(S^{2}(K_1)\times \{r_3\},\rho(r_3)G^*\ame)$
is isometric to a unit sphere.  Thus the curvature of 
$(S^{2}(K_1)\times \{r_3\},G^*\ame)$ is identically $K_3$.  
The case of $n\ge 4$ is
easier.  Once $r_3$ is small enough that 
$(S^{n-1}(K_1)\times \{r_3\},G^*\ame)$ 
contains an open set which in turn contains a closed set isometric to a 
hemisphere
in the standard sphere of constant sectional curvature~$K_3$,
induction and Corollary~\ref {sphererigid2}
implies  
$(S^{n-1}(K_1) \times \{r_3\},G^*\ame)$ 
has constant sectional curvature~$K_3$.
Thus equality holds in (\ref{crap}) as claimed.

But if equality holds in (\ref{crap}), then it follows, using that the 
sectional curvature of $(\am,\ame)$ 
is $\le 1$ and the Gauss curvature equation,
that $\W^2(S_{r_3})\equiv u(R_3)^2\W^2(I)$.  Then 
Theorem~\ref{warprigid2} implies 
$$
G:(\bigcup\{S^{n-1}(K_1)\times \{s\} :0\le s \le   r_3\},\elm)\to (\am,\ame)
$$
is a local isometry.  As the point $y_0=\phi(x_0)$ is in the interior of
the set 
$$
G[\bigcup\{S^{n-1}(K_1)\times \{s\} :0\le s\le r_3\}]\cup\phi[S^{n-1}\times[-r_0,r_0]]
$$
and both $\phi$ and $G$ are local isometries this implies the sectional
curvature of $(\am,\ame)$ is identically $+1$ in a neighborhood of $y_0$.
But this contradicts the choice of $y_0\in\cl(B)$ and completes the proof
 the sectional curvature of $(\am,\ame)$ is identically~$+1$.  This in
turn completes the proof of the theorem.
\end{pf*}

\subsection{Examples}
\label{sec:examples}

We give some examples to show that at least with regard to 
the size of the sets $B$ in the results of Sections~\ref{sec:rigid}
and~\ref{sec:rigid-flat} the results are close to optimal.  However in
light of the many results 
\cite{Beem-Parker:curvature,Dajczer-Nomizu:curvature1,%
Dajczer-Nomizu:curvature2,Graves-Nomizu:curvature,%
Harris:triangle,Kulkarni:curvature,Nomizu:curvature}
to the effect that a manifold with indefinite metric, satisfying a
one sided bound on the
sectional curvature ($=\la R(X,Y)Y,X\ra/(\la X,X\ra\la Y,Y\ra-\la
X,Y\ra^2)$, 
must have constant
sectional curvature, it is worth first giving examples to show that
there are large numbers of semi-Riemannian manifolds that have
one sided curvature bounds in our sense.  

Let $(M_1,g_1)$ and $(M_2,g_2)$ be complete Riemannian manifolds 
and set $M=M_1\times M_2$, $g=g_1-g_2$.  Then $(M,g)$ is a
geodesically complete semi-Riemannian manifold of index $k=\dim M_2$.
If the sectional curvatures of $(M_1,g_1)$ are $\ge0$ and the
sectional curvatures of $(M_2,g_2)$ are $\le 0$ then the curvature of
$(M,g)$ satisfies $R\ge 0$.  (The curvature of $-g$ will satisfy
$R\le0$.) 

For a related example let $(M_1,g_1)$ be a Riemannian
manifold and let $(\h^k,\hm)$ be the hyperbolic space of constant
sectional curvature $-1$.  Let $\rho:\h^k\to[0,\infty)$ be the
Riemannian distance from some point $x_0\in\h^k$.  Let $M=M_1\times
\h^k$ and let $g$ be the warped product metric
$g=\cosh^2(\rho)g_1-\hm$.  If $(M_1,g_1)$ is compact with sectional
curvatures $\ge+1$ then $(M,g)$ is geodesically complete and has
curvature $R\ge 1$.  (The metric $-g$ has curvature $\le -1$.)

To get some more interesting examples we first consider the two
dimensional case. Recall (cf.\ Section~\ref{subsec:two-dim}) that if
$(M,g)$ is two dimensional the Gaussian curvature $K$ is defined by
$K=g(R(X,Y)Y,X)$ where $X$, $Y$ is an orthonormal basis of $T(M)$.
The proof of the following is straightforward.
\begin{prop}\sl\label{prop:2d-example}
Let $\eps_1$, $\eps_2$ be $\pm 1$ and let $g_0=\eps_1dx^2+\eps_2dy^2$
be the standard flat metric on $\r^2$.  Then there are geodesically
complete metrics $g_+$ and $g_-$ on $\r^2$, with Gaussian curvatures
$K_+$ and $K_-$ respectively, so that $g_\pm=g_0$ on the
set $\{y\le 1\}$, $K_+\ge 0$, $K_-\le
0$ and neither of $K_+$ or $K_-$ is identically 0.~\done
\end{prop}
By considering product metrics 
$(\r^{n-2}_j\times \r^2,g_{\r^{n-2}_j}+g_\pm)$ we see that for each
$l$ there are
geodesically complete metrics $g=g_{\r^{n-2}_j}+g_\pm$ on $\r^n$ that
agree with the standard metric
$g^n_k=-\sum_{i=1}^kdx_i^2+\sum_{i=k+1}^ndx_i^2$, on the set $\{x_l\le
1\}$, whose curvature satisfies $R\ge 0$ (or $R\le 0$) but which are
not globally flat.  In this case the exceptional set $B=\{x_l\ge 1\}$
does not satisfy the compact intersection property with respect to any
$\Lambda_\eta$ which shows this condition is necessary in
Theorem~\ref{flatrigid}.  

We note that in the Lorentzian case, the metrics constructed above do
not satisfy the dominant energy condition. This is a consequence of
the following result. Let 
$T_{\Lambda} = \Ric - \half \Scal g +\Lambda g$.  Then the Lorentzian
manifold $(M,g)$ satisfies the {\em dominant energy condition}\/
\cite[\S 4.3]{Hawking-Ellis} if $T_{\Lambda}(X,X) \geq 0$ and 
$Y \mapsto T_{\Lambda}(X,Y)$ defines a nonspacelike
covector for all nonspacelike $X,Y$.

\mnote{check that HE cover $n\geq 4$}
\begin{thm}\sl\label{thm:EinRigid}
Let $(\am,\ame)$ be a timelike geodesically complete globally hyperbolic
Lorentzian manifold of dimension $n$, satisfying the dominant energy condition 
with cosmological constant $\Lambda =  \frac{-(n-1)(n-2)}{2}K_0$ for $K_0 \geq 0$. 
Then, if there is an
isometric immersion of $\r_0^{n-1}(K_0)$ as a totally geodesic hypersurface 
of $\am$, there is a surjective local isometry $\phi: (\r_1^n(K_0)) \to \am$.
\end{thm}
\begin{pf} 
In case $K_0 \geq 0$, $\r_1^n(K_0)$ is globally hyperbolic and
timelike geodesically complete. Hence the result follows from the
conservation theorem (cf.~\cite[p.~94]{Hawking-Ellis}) and uniqueness
theorems for Einstein's equations (cf.~\cite[Ch. 7]{Hawking-Ellis} or
\cite[Ch.~10]{Wald:book}).
\end{pf}

\begin{remark}\normalshape
In case $K_0 >0$, the condition that $(\am,\ame)$ is globally
hyperbolic is not necessary in Theorem~\ref{thm:EinRigid}.
The examples constructed above are globally hyperbolic. 
Therefore Theorem~\ref{thm:EinRigid} implies there is a geodesically 
complete Lorentzian
metric $g$ on $\r^n$ that agrees with the flat metric
$-dx_1^2+dx^2_2+\cdots+dx_n^2$ on the set $\{x_1\le 1\}$, whose
curvature satisfies $R\ge0$, but which is not flat.  Then this metric 
will also satisfy the strong energy condition $\Ric(T,T)\ge 0$ for all
timelike vectors $T$.  However by Theorem~\ref{thm:EinRigid} the
dominant energy condition (with $\Lambda=0$) cannot hold on all of
$(\r^n,g)$.~\done  \end{remark}

Several of the other model spaces we consider are warped products.  To
be concrete consider $\ame=e^{2t}g_0+dt^2$ where $g_0$ is the flat
positive definite metric on $\r^{n-1}$.  Then $\ame$ is the metric on
the hyperbolic space with constant sectional curvature~$-1$.  Choose a
smooth function $w(t)$ on $\r$ so that $w(t)=e^t$ for $t\le 1$ and 
$-\frac{w''(t)}{w(t)}< -1$, and $-\(\frac{w'(t)}{w(t)}\)<-1$ for
$t>1$  (such functions are not hard to construct).  Then the curvature
$\ar_1$ of the metric $\ame_1:=w(t)^2g_0+dt^2$ will satisfy $\ar_1\le -1$ 
but $\ame_1$ is not isometric to the metric $\ame$ even though
they agree on the set $\{t\le 1\}$.  In this case the exceptional set
we are trying to extend across is $B=\{t\ge 1\}$ and this does not have
the compact intersection property with respect to the warping function~$t$.
In this case it is also possible to find functions so that 
$w(t)=e^t$ for $t\le 1$ and $-\frac{w''(t)}{w(t)}> -1$, 
and $-\(\frac{w'(t)}{w(t)}\)>-1$ for $t>1$ which gives examples with
$\ar_1\ge-1$.  Other examples relevant to
Theorem~\ref{thm:rigidity-curved} can be constructed along the same
lines.

\section{Applications}\label{sec:applications}

\subsection{Ends of Constant Curvature and Rigidity}\label{sub:ends}

\newcommand{\mcf}{\mathop{\mathrm{MCTGS}}}     
\newcommand{\bs}{\backslash}            
\newcommand{\f}{\partial}               %
\newcommand{\wh}{\widehat}		

In this section we find all the geodesically complete ends of constant
sectional curvature and with finite fundamental group. This can be
combined 
with our earlier rigidity results to prove rigidity results for
semi-Riemannian manifolds with an end of constant sectional curvature
and finite fundamental group.  

\subsubsection{Structure of Ends of Constant Curvature}

We first remark that in the case of flat ends of complete Riemannian
manifolds there is a structure theory due to Eschenburg and
Schroeder~\cite{Eschenburg-Schroeder:ends} that gives a complete
classification of the locally Euclidean ends.  It would be interesting
to have a corresponding structure theory in the case of ends with
constant sectional curvature~$-1$.  There does not seem to be much
known in this case except when the fundamental group of the end
is finite and then Theorem~\ref{finite-end} below applies.

Let $(\mod, \mme)$ be one of our model spaces and let $G$ be a finite
group of isometries acting on $(\mod,\mme)$.  Assume that there is a
compact set $C\subset \mod$ so that $G$ is fixed point free on
$\mod\setminus C$. (That is if $a\in G$, $a\ne1$ then $a(x)\ne x$ for
all $x\in \mod\setminus C$.  In this case we say that $G$ is {\em fixed
point free on the complement of a compact set\/}.  As the group $G$ is
finite by replacing $C$ with $\bigcup_{a\in G}a[C]$ we can assume that
$C$ is setwise invariant under the action of $G$.  Then the quotient
space $G\bs(\mod\setminus C)$ is a semi-Riemannian manifold with
constant sectional curvature.  We will show that all ends of constant
curvature and finite fundamental group are isometric to examples of
this type and give more information about what groups $G$ can be
realized and how these groups act on $(\mod,\mme)$.
\newline

Let $k$ be an integer, $0 \leq k \leq n$ and let 
$\r^n_k=\r^n_k(0)$ be the geodesically complete flat space form of
index~$k$.  For $k > 0$, write the tangent space to $\r_k^n$ as 
$T(\r^n_k)_0=\r^k_k\oplus\r^{n-k}_0$ and let
$S^{k-1}\times S^{n-k-1}$ be the subset of
$T(\r^n_k)_0=\r^k_k\oplus\r^{n-k}_0$ defined by 
$\{(x,y):\la x,x\ra=-1, \la y,y\ra=1\}$.  
We view $S^{k-1}\times
S^{n-k-1}$ as the Riemannian product of spheres. 

Let $O(k,n-k)$ denote the group of isometries of $\r^n_k$ fixing the origin. 
This can be identified with the  group of linear isometries of the underlying 
indefinite inner product space.
Note  that any Riemannian isometry $a$ of $S^{k-1}\times S^{n-k-1}$ which 
preserves the product structure (i.e. $a\in O(k)\times O(n-k)$) is induced by a
unique $g \in O(k,n-k)$,   
and the action of the derivative $g_*$ on the
tangent space $\r_k^n\oplus\r^{n-k}_0$ preserves the splitting and
induces the mapping $a$ on $S^{k-1}\times S^{n-k-1}$.  Thus any subgroup
$G$ of $O(k)\times O(n-k)$ extends to 
to a subgroup of $O(k,n-k)$.

Call a
subgroup $G$ of $O(k)\times O(n-k)$ {\em strongly  fixed point free\/}
iff for every $a=(a_1,a_2)\in G$ if there is a point $u_1\in S^{k}$
with $a_1u_1=u_1$, or a point $u_2\in S^{n-k}$ with $a_2u_2=u_2$ then
both $a_1$ and $a_2$ are the identity.  That is $G$ is strongly fixed
point free iff both the induced actions on $S^k$ and $S^{n-k}$ are
fixed point free.

It is easy to check that if $G$ is a subgroup of
$O(k)\times O(n-k)$ then the induced action of $G$ on $\r^n_k$ is
fixed point free on $\r^n_k\setminus \{0\}$ iff the action is strongly
fixed point free.  If $G$ is finite and strongly fixed point free, then 
we can 
form the orbit space $G\bs \r_k^n$.  This  is a smooth manifold except
at the orbit corresponding to the origin which is a conical singular point.

\begin{definition}\sl Let $G\subset O(k)\times O(n-k)$ be a finite 
group of isometries of the
space $S^{k-1}\times S^{n-k-1}$ which is strongly fixed point free.
Then the 
{\em standard end\/} of $\r^n_k$ determined by $G$ is the end of the
space $G\bs \r^n_k$ constructed above.~\done
\end{definition}

Let $K_0>0$.  To make notation easier normalize so that $K_0=1$. 
Let $(\h^k,\hm)$ be the hyperbolic space with constant sectional
curvature~$-1$ and let $(S^{n-k},\elm)$ be the standard sphere.
Choose a point $x_0\in\h^k$ to use as an origin and let $\rho:\h^k\to
[0,\infty)$ be the Riemannian distance from $x_0$.  

For  $1 \leq k <  n-1$, the space
$\r^n_k(+1)$ has a representation as the warped product 
$(S^{n-k}\times\h^k,\cosh^2(\rho)\elm-\hm)$.  Let $S^{k-1}$ be the
unit sphere of $T(\h^k)_{x_0}$.  Any isometry $a$ of $S^{k-1}$ is
induced by a unique isometry of $\h^k$ fixing the origin $x_0$.  Thus
any isometry $a=(a_1,a_2)\in O(n-k+1)\times O(k)$ of $S^{n-k}\times
S^{k-1}$ induces an isometry of $\r^n_k(+1)=
(S^{n-k}\times\h^k,\cosh^2(\rho)\elm-\hm)$ in the obvious way.  

If $G\subset O(n-k+1)\times O(k)$ we say  $G$ is {\em weakly fixed
point free\/} iff for any non-identity 
$(a_1,a_2)\in O(n-k+1)\times O(k)$ if $a_2$ fixes a point of
$S^{k-1}$, then $a_1$ has no fixed points on $S^{n-k}$.
Note that  if $G\subset O(n-k+1)\times O(k)$ 
is a finite group, then the induced action on
$\r^n_k(+1)$ is fixed point free on the complement of some compact set
iff $G$ is weakly fixed point free.  
Thus the quotient space $G\bs
\r^n_k(+1)$ is smooth away form a compact set.

Finally consider the case of $k=n-1$, and let $x_0$ be a point of
$\r^n_{n-1}(K_0)$ with $K_0>0$.  Let $\Bbb Z_2$ act on $\r^n_{n-1}(K_0)$ by the geodesic symmetry through $x_0$.
In this case the only standard ends of $\r^n_{n-1}$ are the end of
$\r^n_{n-1}$ and the end of the orbit space $\Bbb Z_2\bs
\r^n_{n-1}(K_0)$.

\begin{definition}\sl Let $K_0>0$ and $1 \leq k <  n-1$
Let $G\subset O(n-k+1)\times O(k)$ be a weakly fixed point free finite
group.  Then the end of the quotient $G\bs \r^n_k$ constructed above
is a {\em standard end\/}. If $k=n-1$, then the standard end is as just
described.  If $K_0<0$ then the space
$\r^n_k(K_0)$ is anti--isometric to the space $\r^n_{n-k}(-K_0)$.  Thus
we use the constructions just given to define a standard end in this 
case.~\done
\end{definition}

Let $E$ be an end of the manifold $M$.  Recall that the fundamental
group $\pi_1(E)$ is isomorphic to $G$ iff for every compact set
$C_1\subset M$ there is a compact set $C_2\supseteq C_1$ so that
$\pi_1(E\setminus C_2$ is isomorphic to $G$ and if $i:(E\setminus
C_2)\to (E\setminus C_1)$ is inclusion, then the induced map
$i_*:\pi_1(E\setminus C_2)\to \pi_1(E\setminus C_1)$ is injective.  In
particular this means  the end $E$ is simply connected iff of all
compact $C_1$ there is a $C_2\supseteq C_1$ so that $E\setminus C_2$
is simply connected.  

In some cases the model space $\r^n_k(K_0)$ will have two ends.  In
this case the two ends are isometric
and so the isomorphism class of
an end of $\r^n_k(K_0)$ is well defined.  

A more troublesome case is when the end of the model space is not simply
connected.  (This happens when $\r^n_k(K_0)$ is diffeomorphic to
$\r^2$ or 
$S^{n-2}\times \r^2$ where $n\ge 4$.)
These facts and some other
information we will need is given in the following table.  The column
headed by ``Diff. Type'' gives the diffeomorphism type of
$\r^n_k(K_0)$, the column headed by ``Positive'' (resp.\ ``Negative'')
gives the lengths of the positive (resp.\ negative) geodesics (when
this number is finite it means that all geodesics are closed with the
given length as smallest period).  The
other two columns give the number of ends and the fundamental group
of an end.

\begin{table}[htp]
\caption[]{}
\label{tablepi}
\begin{center}
\begin{tabular}{|c|c|c|c|c|c|c|}\hline
$K_0$&$k$&Diff.\ Type&Positive&Negative& \# of ends&$\pi_1(\infty)$\\
\hline\hline
$K_0=0$&$1\le k\le n-1$&$\r^n$&$\infty$&$\infty$&1&$\la0\ra$\\
\hline
$K_0=0$&$ k=0$&$\r^n$&$\infty$&---&1&$\la0\ra$\\
\hline
$K_0=0$&$ k=n $&$\r^n$&---&$\infty$&1&$\la0\ra$\\
\hline
$K_0<0$&$ k=0$&$\r^n$&$\infty$&---& 1&$\la 0\ra$\\ 
\hline
$K_0<0$&$ k=1$&$\r^n$&$\infty$&$\infty$& 1&$\la 0\ra$\\ 
\hline
$K_0<0$&$2\le k=n-1$&$S^{n-1}\times\r$&$\infty$&$2\pi/\sqrt{K_0}$& 2
	&$\la0\ra$\\
\hline
$K_0<0$&$2\le k\le n-3$&$S^{k}\times\r^{n-k}$&$\infty$&$2\pi/\sqrt{K_0}$
						& 1&$\la0\ra$\\
\hline
$K_0<0$&$2\le k= n-2$&$S^{k}\times\r^{2}$&$\infty$&$2\pi/\sqrt{K_0}$
						& 1&$\Bbb Z$\\
\hline$K_0<0$&$k=n$&$S^{n}$&---&$2\pi/\sqrt{K_0}$
						& 0&---\\
\hline
$K_0>0$&$k=0$&$S^n$&$2\pi/\sqrt{K_0}$&---&0&---\\ 
\hline
$K_0>0$&$k = n$&$\r^n$&---&$\infty$&1&$\la0\ra$\\
\hline
$K_0>0$&$ k= n-1$&$\r^n$&$\infty$&$\infty$&1&$\la0\ra$\\
\hline
$K_0>0$&$k=1,n\ge3$&$S^{n-1}\times\r$&$2\pi/\sqrt{K_0}$&$\infty$
		&2&$\la0\ra$\\
\hline
$K_0>0$&$3\le k\le n-2$&$S^{n-k}\times\r^k$&$2\pi/\sqrt{K_0}$&$\infty$
		&1		&$\la0\ra$\\
\hline
$K_0>0$&$2= k\le n-2$&$S^{n-k}\times\r^2$&$2\pi/\sqrt{K_0}$&$\infty$
		&1		&$\Bbb Z$\\
\hline
\end{tabular}
\end{center}
\end{table}

We will use the fact that when a geodesic has infinite length then it
will eventually leave every compact set.  Null geodesics have this
property in all cases.  

\begin{thm}\sl \label{finite-end}
Let $(M,g)$ be a geodesically complete semi-Riemannian
manifold of dimension~$n$ and index~$k$ and let $E$ be an end of $M$.
Assume that $E$ has finite fundamental group and constant sectional
curvature $K_0$.  Then the model  space $\r^n_k(K_0)$ is non-compact,
its ends are simply connected, and the end $E$ is isometric  to a
standard end of constant sectional curvature~$K_0$ and index~$k$.
\end{thm}

This implies that when the model space $(\mod,\mme)$ is diffeomorphic
to $\r^2$ or $S^{n-2}\times \r^2$ for $n \geq 4$
(the cases where the end of $\mod$
has infinite fundamental group) then no end with a finite
fundamental group can be geodesically complete and be locally
isometric to an end of $(\mod,\mme)$.
\newline

The proof of Theorem~\ref{finite-end}
breaks  up into two parts.  The first is understanding the
simply connected ends of constant curvature and the second is in
proving a fixed point theorem for finite groups of isometries on the
model space.  The facts  about simply connected ends are given  by

\begin{prop}\sl \label{simple-end}
Let $(M,g)$  be a geodesically complete semi-Riemannian
manifold of dimension $n$ and index $k$. Assume that $M$ has non--empty
compact boundary $\partial M$ and that $(M,g)$ is geodesically complete in
the sense that every inextendible geodesic either has no endpoint or an
endpoint on $\partial M$. 

Assume  that $(M,g)$ is simply
connected and that it has constant sectional curvature~$K_0$ and
that the model space $\mod$ has simply connected ends. 
Then $k$ and $K_0$ are such that $\mod$ is noncompact and
there are compact sets $\f M\subset C_1\subset M$  and
$C_2\subset\r^n_k(K_0)$ so that every unbounded component of
$M\setminus C_1$ is isometric with an unbounded component of
$\r^n_k(K_0)\setminus C_2$.
\end{prop}

\begin{pf} As $M$ is
simply connected and locally isomorphic to $\r^n_k(K_0)$ a standard
monodromy argument shows there is a local isometry $f:M\to
\r^n_k(K_0)$.
Let $U$ be a connected component of $\r_k^n(K_0)\setminus f[\f M]$ and
set $\wh{U}:=f^{-1}[U]$.  If $\wh{U}\ne \emptyset$ then
$f|_{\wh{U}}:\wh{U}\to U$ is a covering map.  To see this let $y\in U$
and let $B\subset U $ be a geodesically convex neighborhood of $y$
whose closure is disjoint from $f[\f M]$ and let $N_y\subset T(\r_k^n(K_0))_y$
be the convex neighborhood of the origin so that $\exp[N_y]=B_y$.  Let 
$x\in \wh{U}$ so that $f(x)=y$.   Let $\wh{N}_x\subset T(M)_x$ be the
convex neighborhood of the origin so that $f_{*x}[\wh{N}_x]=N_y$ and let 
$\wh{B}_x=\exp[\wh{N}_x]$.  Then $\wh{B}_x\subset \wh{U}$, and
$f|_{\wh{B}_x}:\wh{B}_x\to B_x$ is a diffeomorphism.  As $x$ was any
point of $\wh{U}$ with $f(x)=y$ this shows that the neighborhood $B_x$ is
evenly covered.  Thus if $\wh{U}\ne\emptyset$ the map $f|_{\wh{U}}$ is
a covering map as claimed.
\smallskip

\noindent
{\bf Claim:\ }  The set $f[M]$ is not contained in an compact
subset of $\r_k^n(K_0)$.

Assuming the claim we prove the proposition.  
Let $U$ be an unbounded component of $\r^n_k(K_0) \setminus f[\partial M]$
so that $f[M]\cap U\ne\emptyset$ and set $\wh U:=f^{-1}[U]$.  
Then $f|_{\wh U}:\wh U\to
U$ is a covering map.  
Then there is a
compact set $C_2\supset f[\f M]$ so that $U_2:=U\setminus C_2$ is
simply connected.  Let $\wh{U}_2$ be a connected component of
$f^{-1}[U_2]$.  As $U_2$ is simply connected the map
$f|_{\widehat{U}_2}:\widehat{U}_2\to U_2$  must be a
diffeomorphism and thus an isometry.  As $\f U_2$ is compact the
boundary of $\wh{U}_2$ is also compact.  The proposition now follows.
\smallskip

We now prove the claim.  
This is complicated by the fact that in many cases the
model space $\r^n_k(K_0)$ contains closed geodesics.  This means it
might be possible to have a geodesic ray $\gamma:[0,\infty)\to M$ that
eventually leaves every compact set of $M$, but $f\circ\gamma:[0,\infty)\to
\r^n_k(K_0)$ is closed and thus is contained in a compact subset of
$\r^n_k(K_0)$.  To avoid this note that by multiplying the metric of
$\r^n_k(K_0)$ by $\pm 1$  we can assume that no non-negative geodesic
of $\r^n_k(K_0)$ is closed and that in fact a non-negative geodesic eventually
leaves every compact set.

Let $D$ be a compact subdomain of $M$
with smooth boundary and so that $\f M$
is a subset of the interior of $D$.  We also choose a complete
Riemannian metric $h$ on $M$ (which has no relation to the semi-Riemannian
metric $g$ on $M$).  Let $S^+(D)$ be the set of all non-negative (i.e.
$g(u,u)\ge 0$) $h$-unit vectors tangent to $M$ at a point of $D$.
Then all the fibers of $S^+(D)$ are compact and $D$ is compact so
the space $S^+(D)$ is compact.  Let $V$ be an unbounded
component 
of $M\setminus D$ and let $\cal G$ be the compact set of ordered
pairs $(u,t)\in S^+(D)\times (0,\infty)$ so that $\exp(tu)\in V$ and
$0\le t\le1$.
(Here $\exp$ is the exponential of the semi-Riemannian metric $g$.)
Note that if $x\in \f D\cap \f V$ then there is a positive $h$-unit
vector $u$ at $x$ pointing into $V$
and thus $\exp(tu)\in V$ for small
$t$ and thus $\cal G$ is not empty.  For each $(u,t)\in {\cal G}$  let
$(\ga(u,t),\gb(u,t))$ be the maximal interval containing $t$ and
contained in the set $\{s:\exp(su)\in V\}$.  Thus $0\le
\ga(u,t)<t<\gb(u,t)\le \infty$ and $\ga(u,t)\le 1$.  
If for some $(u,t)\in {\cal G}$ we have
$\gb(u,t)=\infty$ then $\gamma(s)=f(\exp(su))$ with $s>\ga(u,t)$ 
is a non-negative geodesic ray of $\r^n_k(K_0)$ and thus it eventually
leaves every compact set.  Thus $f[M]$ is not contained in any compact
subset of $\r^n_k(K_0)$ and thus the claim holds in this case.

Next assume that there is a sequence
$\{(u_l,t_l)\}_{l=1}^\infty\subset {\cal G}$ so that
$\lim_{l\to\infty}\gb(u_l,t_l)=\infty$.  Then as $t_l\in[0,1]$ and $S^+(D)$ is
compact, by going to a subsequence we can assume  $t_l\to t_0$ for
some $t_0\in[0,1]$ and $u_l\to u_0$ for some $u_0\in S^+(D)$.  Consider
the geodesic $\gamma(s)=\exp(su_0)=\lim_{l\to\infty}\exp(su_l)$ for
$s>1$.  This will be disjoint from $D^\circ$ (the interior of $D$) and
thus does not meet $\f M$.  It is also a non-negative geodesic.  
Thus the geodesic $f(\gamma(s))$ in
$\r^n_k(K_0)$ will eventually leave every compact set and thus $f[M]$ is
not contained in any compact subset of $\r^n_k(K_0)$.  

This leaves the case where the function $\gb(u,t)$ is bounded on $\cal
G$, say $\gb(u,t)\le C_0$.  Then let 
$$
G=\{\exp(su): \mbox{$\ga(u,t)<s<\gb(u,t)$ for some $(u,t)\in {\cal G}$}\}.
$$
Consider the sequence $\{\exp(s_lu_l)\}_{l=1}^\infty\subset G$ 
where $(u_l,t_l)$ is a sequence in  $\cal G$.
By passing to a subsequence we can
assume  that $t_l\to t_0$ and $u_l\to u_0$ for some $u_0\in S^+(D)$.
As $s_l<\gb(u_l,t_l)\le C_0$ we can also assume  $s_l\to s_0$ for
some $s_0\le C_0$.   Then $\exp(s_lu_l)\to \exp(s_0u_0)$ which shows
that every sequence out of $G$ has a convergent subsequence.
Therefore so the 
closure of $G$ is compact in $M$.  As the component $V$ is unbounded it
can not be contained in any compact set and thus there is a  point $x$
in $V$ that is not in $G$.  Let $\gamma$ be a positive geodesic of
$M$ through the point $x$.  Then $\gamma$ can not meet $D$, for if it
did the point $x$ would be of the form $\exp(su)$ for some
$(u,t)\in{\cal G}$  contradicting that $x$ is not in $G$.  Thus
the geodesic $f\circ\gamma$ of $\r^n_k(K_0)$ leaves every compact set
and completes the proof of the claim $f[M]$ is not contained in any 
compact subset of $\r^n_k(K_0)$.  This also completes the proof of the 
proposition.
\end{pf}

To pass from the case of simply connected ends to the case of ends
with finite fundamental group we need information about finite group
actions on certain homogeneous spaces of totally geodesic
submanifolds of $\r^n_k(K_0)$.  By a {\em maximal compact totally geodesic
submanifold\/} in
$\r^n_k(K_0)$ we mean a compact totally geodesic submanifold $C$ of
$\r^n_k(K_0)$ which is not properly contained in any other compact
totally geodesic submanifold of $\r^n_k(K_0)$.  In many cases the
maximal compact totally geodesic submanifolds may just be points.  
In the case the space
$\r^n_k(K_0)$ is compact, then the only maximal compact flat is
$\r^n_k(K_0)$ itself.  The following gives the basic facts.  The
proof is left to the reader.

\begin{prop}\sl If $C_1$ and $C_2$ are maximal compact totally 
geodesic submanifolds in
$\r^n_k(K_0)$, then there is an isometry $g$ of $\r^n_k(K_0)$ with
$gC_1=C_2$.  Thus all the maximal compact totally geodesic submanifolds in $\r^n_k(K_0)$ are
isometric.  The maximal compact totally geodesic submanifolds are
\begin{table}[htp]
\caption[]{}
\label{table3}
\begin{center}
\begin{tabular}{|c|c|c|}\hline
$K_0$&$n,k$& Max. Cpt. Tot. Geo. Subman.\\
\hline\hline
$K_0=0$&$0\le k\le n$& Point\\ \hline
$K_0>0$&$k=n-1,n$&Point\\ \hline
$K_0>0$&$0\le k\le n-2$& $S^{n-k}(1/\sqrt{K_0})$\\ \hline
$K_0<0$&$k=0,1$& Point\\ \hline
$K_0<0$&$2\le k\le n$& $S^k(1/\sqrt{-K_0})$\\\hline
\end{tabular}
\end{center}
\end{table}
\end{prop}

For the rest of this section we will denote the space of maximal
compact totally geodesic submanifolds in $\r^n_k(K_0)$ by $\mcf(n,k,K_0)$.
Note that if the elements of $\mcf(n,k,K_0)$ are points, then 
$\mcf(n,k,K_0) = \mod$.

\begin{prop}\sl $K_0>0$ and $k\ne n-1$ or $K_0<0$ and $k\ne 1$ then 
$\mcf(n,k,K_0)$ can be given the structure of a Riemannian symmetric
space of non-compact type in such a way that the group of isometries of
$\r^n_k(K_0)$ acts on $\mcf(n,k,K_0)$ by isometries.  In particular the
sectional curvature on $\mcf(n,k,K_0)$ is non-positive.
\end{prop}

\begin{pf}  As the spaces $\r^n_k(K_0)$ and $\r^n_{n-k}(-K_0)$ are
anti--isomorphic, it is enough to prove the result in the case
$K_0>0$.  When $K_0>0$ the full isometry group of $\r^n_k(K_0)$ is the
orthogonal group $O(n-k+1,k)$.  If $C_0\in \mcf$ then the subgroup of
$O(n-k+1,k)$ fixing $C_0$ is $O(n-k+1)\times O(k)$ and thus as a
homogeneous space $\mcf=O(n-k+1,k)/(O(n-k+1)\times O(k))$.  This is
well known to be a Riemannian symmetric space
(\cite[p518]{Helgason:DS} it is the
non-compact dual to the Grassmannian manifold
$G_k(\r^{n+1})=O(n)/(O(n-k+1)\times O(k))$). It is
a standard result that $\mcf=O(n-k+1,k)/(O(n-k+1)\times O(k))$ has
non-positive curvature as it is a Riemannian symmetric space of
non-compact type.  Cf.~\cite[Theorem~3.1~p. 241]{Helgason:DS}.
\end{pf}

\begin{prop}\sl \label{fixed-flat}
Let $H$ be a compact group acting on $\r^n_k(K_0)$ by
isometries.  Then $H$ fixes some point of $\mcf(n,k,K_0)$.  
\end{prop}
\begin{pf}
If $K_0=0$ then $H$ is a compact group acting 
on $\r^n=\r^n_k(0)$ by
affine maps.  Thus $H$ must fix a point of $\r^n_k(K_0)$ which proves
the result in this case as the maximal compact totally geodesic submanifolds 
are points.

If $K_0>0$ and $k\ne n-1$ of $K_0<0$ and $k\ne 1$ then by the last
proposition $\mcf(n,k,K_0)$ is a simply connected Riemannian manifold
with non-positive sectional curvature and $H$ acts on $\mcf(n,k,K_0)$
by isometries.  So by the Cartan fixed point theorem \cite[Theorem
13.5, p. 75]{Helgason:DS} $H$ fixes a point
of $\mcf(n,k,K_0)$.

This leaves the equivalent cases of $K_0>0$ and $k=n-1$ and $K_0<0$
and $k=1$.  We will deal with the case where $K_0<0$ and $k=1$ and we
normalize so that $K_0=-1$.  Let $\r^{n+1}_2$ be $\r^{n+1}$ with the inner
product $\la x,y\ra=-x_1y_1-x_2y_2+x_3y_3+\cdots+x_{n+1}y_{n+1}$.  
Let $S_1^{n}(-1)$ be the hypersurface in $\r^{n+1}_2$ defined by $\la
x,x\ra=-1$. Then $S_1^n(-1)$ has constant sectional curvature $-1$ and
is diffeomorphic to $S^1\times \r^{n-1}$. 

Let $p:\r^n_1(-1)\to
S_1^n(-1)$ be the covering map and let $\Gamma$ be the group of deck
transformations.  
Recall that  $\r_1^n(-1)$ is time orientable.  
Fix a time orientation on $\r^n_1(-1)$.  Because each timelike
geodesic of $S_1^n(-1)$ is closed of length~$2\pi$ it is not hard to see
that the group $\Gamma$ of deck transformations is cyclic and its
generator is the unique isometry $a$ of $\r^n_1(-1)$ with the property
that for every future pointing timelike unit speed geodesic $\gamma$
in $\r^n_1(-1)$ there holds $a\gamma(t)=\gamma(t+2\pi)$.  

Let $g$ be
any isometry of $\r^n_1(-1)$ that preserves the time orientation of
$\r_1^n(-1)$.  Then for any future pointing time like unit speed 
geodesic $\gamma$ the curve $g\gamma$ is also a future pointing time 
like unit speed geodesic and so
$ag\gamma(t)=g\gamma(t+2\pi)=ga\gamma(t)$.  As there is a 
future pointing time like unit speed  geodesic through any point of
$\r^n_1(-1)$ this implies $ag=ga$.  That is the generator $a$ of
$\Gamma$ commutes with every isometry of $\r_1^n(-1)$ that preserves
the time orientation.  

We now prove the proposition in the special case
where every element of the group $H$ preserves the time
orientation of $\r^n_1$.  Then every element of $H$ commutes with every
element of the group of deck transformations and therefore $H$ also
acts on $S_1^n(-1)$ by $gp(x):=p(gx)$.  
We now claim that there is a time like geodesic of $S_1^n(-1)$
that is fixed by the action of $H$.  Note that the group of isometries
of $S_1^n(-1)$ that stabilize some fixed timelike geodesic is
$O(2)\times O(n-1)$ and the full isometry group of $S_1^n(-1)$ is
$O(2,n-1)$.
Therefore as a homogeneous space, the space of timelike
geodesics of $S_1^n(-1)$ is $O(2,n-1)/(O(2)\times O(n-1))$ which, just
as above, is a Riemannian symmetric space with non-positive sectional
curvature.  Thus the Cartan fixed point theorem lets us conclude 
$H$ fixes a timelike geodesic $c$ of $S_1^n(-1)$.  

Let $\gamma$ be the
future pointing unit speed time like geodesic covering of $S_1^n(-1)$
covering $c$.  Then $\gamma$ is fixed set-wise by every element $g\in
H$.  As the elements of $H$ preserve the time orientation this implies
that if $g\in H$, then $g\gamma(t)=\gamma(t+t_0)$ for some
$t_0\in\r^n$.  But then $g^l\gamma(t)=\gamma(t+lt_0)$.  But as the
group $H$ is compact the orbit $\{g^l\gamma(0)=\gamma(lt_0): l=0,
\pm1,\pm2,\ldots\}$ must have compact closure in $\r^n_1(-1)$ which is
only the case if $t_0=0$.  Thus $g\gamma(t)=\gamma(t)$ for all $g\in
H$.  So in the case  where all elements of $H$ preserve the time
orientation of $\r^n_1(-1)$ there is always a time like geodesic that
is fixed pointwise by $H$.

Finally we consider the case where the group contains elements that
reverse the time orientation of $\r^n_1(-1)$.  In this case let $H_0$ be
the subgroup of $H$ of elements that preserve the the time orientation
of $\r^n_1(-1)$.  Then  $H_0$ is a subgroup of $H$ of index~2.  We
consider two cases.  First assume  every element of $H_0$ fixes
every element of $\r_1^n(-1)$.  Then $H_0=\la1\ra$ and so $H$ is just
a cyclic group of order two.   Let $b$ be the generator of $H$.  Then
$b$ reverses the time orientation of $\r_1^n(-1)$.  Choose a fixed
future pointing time like unit speed geodesic $\gamma_0$ for
reference and let $\tau$ be the time function
determined by $\gamma_0$,  
cf.~Definition \ref{def:time-fcn}.
That is, on $\gamma_0$ the function is given by $\tau(\gamma_0(t))=t$ 
and then
$\tau$ is extended to $\r_1^n(-1)$ to be constant on the totally
geodesic hypersurfaces perpendicular to $\gamma_0$.  As the element
$b$ of $H$ reverses the time orientation of $\r^n_{1}(-1)$ the
function $\tau\circ b$ is a decreasing function when restricted to
$\gamma_0$, $\lim_{t\to\infty}\tau(b\gamma(t))=-\infty$, and
$\lim_{t\to-\infty}\tau(b\gamma(t))=+\infty$.  Thus there is a point
$x=\gamma(t_1)$  so that $\tau(x)=\tau(gx)$.  

Let $M$ be the set of
all points of $\r^n_1(-1)$ where $\tau(x)=\tau(bx)$.  Then $M$ is not
empty and as the level sets of both $\tau$ and $\tau\circ b$ are
totally geodesic spacelike hypersurfaces and $M$ is the intersection of
two such level sets, $M$ will be a space like totally geodesic 
submanifold of either codimension 1 (in the case the level sets 
happen to be equal) or codimension 2 (in all other cases).  This
implies  $M$ is a simply connected Riemannian manifold of
constant non-positive sectional curvature.  As the isometry $b$ leaves
$M$ fixed set-wise, the Cartan fixed point theorem can again be used to
show  $b$ has a fixed point on $M$.  This completes the proof in
the case $H_0=\la1\ra$.

This leaves the case where $H_0\ne\la1\ra$.  Then by what we have shown
above there is a timelike geodesic $\gamma$ that is fixed pointwise by
$H_0$.  Let $N$ be the subset of $\r^n_1(-1)$ of all points fixed by
all elements of $H_0$.  As $N$ is the fixed point set of a group of
isometries it is a totally geodesic submanifold of $\r_1^n(-1)$. As $N$
contains a timelike geodesic the restriction of the metric to $N$ is
non-degenerate (in Lorentzian manifolds the only degenerate subspaces
have a positive semidefinite metric).  As $H_0\ne\la1\ra$ the
submanifold $N$ has dimension less than $n$.
As the subgroup
$H_0$ of $H$ has index two it is normal in $H$.  This implies that the
set $N$ is fixed set-wise by the group $H$.  To see this let $a\in H_0$
and $g\in H$.  Then $g^{-1}ag\in H_0$ and so for $x\in N$ we have 
$agx=gg^{-1}agx=gx$.  Thus $gx$ is fixed by all elements of $H_0$ and so
$gx\in N$ as claimed.  Thus $H$ acts on the submanifold $N$ and $N$ is
isometric to $\r^m_1(-1)$ for some $m<n$.  Thus we can now use
induction to conclude that $H$ has a fixed point in $N$ and therefore
also in $\r_1^n(-1)$.  This completes the proof.
\end{pf}

\begin{pf*}{Proof of Theorem~\ref{finite-end}}  
Let $E$ be an end of $M$
with constant sectional curvature and finite fundamental group
$\pi_1(E)$.  Then there is a compact set $C$, which we can assume to
have smooth boundary, so that $E\setminus C$ has fundamental group
isomorphic to $\pi_1(E)$.  Let $\pi:N\to E\setminus C$ be the
universal cover of $E\setminus C$ and let $G$ be the group of deck
transformations of this cover.  Thus $G$ acts on the semi-Riemannian
manifold $N$ by isometries and $G$ and the orbit space $G\bs N$ is
isometric to $E\setminus C$. Thus it is enough to show  $G\bs N$
is isometric to a standard end.  By use of Proposition~\ref{simple-end} by
replacing $C$ by a larger compact set we can assume  $N$ is
isometric to an unbounded component of $\r^n_k(K_0)\setminus C_2$  for
some compact $C_2$. Thus we simply assume $N=\r^n_k(K_0)\backslash C_2$.
As the space $\r^n_k(K_0)$ is simply connected, the
action of $G$ on $N$ extends to an action on the whole space
$\r^n_k(K_0)$.  By Proposition~\ref{fixed-flat} the group fixes some
element of $\mcf(n,k,K_0)$.  
If the maximal compact totally geodesic submanifolds of
$\r_k^n(K_0)$ are spheres $S^p$ with $p\ge 2$, then by normalizing
the metric we can assume the space is $(S^p\times
\h^{n-p},\cosh^2(\rho)\elm-\hm)$ and that the fixed element of
$\mcf(n,k,K_0)$ is $S^p\times\{\rho=0\}$.   
Then it follows that $G$
must be fixed weakly fixed point free and that $G\bs N$ is a standard
end.

If the maximal compact totally geodesic submanifolds are points, then let 
$x_0$ be the point of $\r^n_k(K_0)$ fixed by $G$.  
Then there is a $k$-dimensional
subspace $V$ of $T(\r^n_k(K_0))_{x_0}$ that is fixed by $G$. (I.e.
$a_{*x_0}V=V$ for all $a\in G$.)   This is because the set of all such
subspaces is just the homogeneous space $O(k,n-k) /(O(k)\times O(n-k)$
and as we have done several times above, we can apply the 
Cartan fixed point theorem to this case.  
But then by standard linear algebra the
orthogonal complement $V\oc$ of $V$ is also invariant under $G$.
Letting $S^{k-1}\times S^{n-k-1}$ be the product of the unit spheres
in these spaces we see that the group has an action on 
$S^{k-1}\times S^{n-k-1}$ which must then be strongly fixed point free
and thus the space $G\bs N$ is again isomorphic to a standard end.
This completes the proof.
\end{pf*}

\subsubsection{Rigidity Results}

\begin{thm}\sl \label{rigid-flat-end} 
Let $(\am,\ame)$ be a geodesically complete semi-Riemannian manifold of
dimension $n\ge3$ and index~$k$ with
curvature satisfying one of the two inequalities $\ar \geq 0$
or $\ar\le 0$  and assume  $(\am,\ame)$ has an end $E$ with
$\ar\equiv 0$ on $E$ and $\pi_1(E)$ finite.  Then $(\am,\ame)$ is isometric
to the flat model space $(\r^n_k,g_0)$.
\end{thm}

\begin{thm}\sl \label{rigid-lorentz-end}
Let $(\am,\ame)$ be a geodesically complete Lorentzian or Riemannian 
manifold of dimension $n\ge3$ with
curvature tensor satisfying one of the two one sided inequalities
$\ar\ge -1$ or 
$\ar\le -1$.  Assume  $(\am,\ame)$ has
an end $E$ with $\pi_1(E)$ finite and so that 
$\ar\equiv -1$ on $E$.  Then if  $(\am,\ame)$ is Lorentzian it is 
isometric to the anti-de~Sitter space $(\r^n_1(-1),g_{\r^n_1(-1)})$. If
$(\am,\ame)$ is Riemannian it is isometric to the hyperbolic space
$(\r^n_0(-1),g_{\r^n_0(-1)})$
\end{thm}

\begin{pf*}{Proofs}  We first consider
Theorem~\ref{rigid-flat-end}.  By Theorem~\ref{finite-end} there is a
compact set $C\subset\r^n_k$ and an local isometry
$\phi:\r^n_k\setminus C\to \am$ so that for each $y\in \am$ the preimage
$\phi^{-1}[y]$ has at most $\#(\pi_1(E))$ elements.  By
Theorem~\ref{flatrigid} this extends to a surjective local isometry
$\widehat{\phi}:\r^n_k\to \am$.  This implies the map $\widehat{\phi}$
is a covering map and that the group of deck transformations is
isomorphic to $G:=\pi_1(E)$ and thus $G$ is finite. By 
Proposition~\ref{fixed-flat} $G$ has a fixed point on $\r^n_k$.
But as $G$ is the group of deck transformations of
$\widehat{\phi}:\r^n_k\to M$ the only way it can
have a fixed point is if it is the trivial group $G=\la1\ra$.  This
implies  $\widehat{\phi}$ is an isometry and completes the proof of
Theorem~\ref{rigid-flat-end}.  The proof of
Theorem~\ref{rigid-lorentz-end} is identical except that the rigidity
results of Section~\ref{sec:rigid}
are used.
\end{pf*}

\begin{remark}\normalshape
In all cases not covered by the last two theorems
the exact analog of the theorems is false.  First consider the case
where $K_0<0$ and $3\le k\le n-1$.  Normalize so that $K_0=-1$.  Then
as above $\r^n_k(K_0)$ is isometric to the warped product metric
$g_0=-\cosh^2(\rho)\elm +\hm$ on $S^k\times\h^{n-k}$. Thus if
$G$ is a finite group of isometries of $(S^k,\elm)$ acting without
fixed points, then it also has a fixed point free action on
$\r^n_k(K_0)\approx S^k\times\h^{n-k}$ by letting it act on the first
factor.   Then the manifold $G\bs \r^n_k$ has a metric of constant
sectional curvature~$-1$, and an end with a finite fundamental group
without being isometric to $\r^n_k(K_0)$.  There are similar examples
when $K_0>0$.

However it is still reasonable to conjecture that if $(M,g)$ is a
geodesically complete semi-Riemannian manifold of dimension~$n\ge3$
and index~$k$ with $3\le k\le n-1$, which satisfies an appropriate one
sided bound on sectional curvature and has an end $G$ with finite
fundamental group and constant sectional curvature $K_0<0$ that
$(M,g)$ is isometric to a quotient of $\r^n_k(K_0)$ by a finite fixed
point free group of isometries.  Getting control of the geometry of
the end is taken care of by Theorem~\ref{finite-end}.  What is needed
to prove this is an extension of the rigidity results of
Section~\ref{sec:rigid} to the case of space forms of index other than
$k=0,1, n-1,n$.~\done
\end{remark}

\subsection{Rigidity in Quotients}\label{sub:quot}

\begin{definition}\sl Let $M$ be a smooth manifold and let $C\subset M$ be
a closed subset.  Let $g_1$ and $g_2$ be two smooth semi-Riemannian metrics on 
$M$ and let $\ell$ be a non-negative integer.  Then we say that 
$g_1$ and $g_2$ agree to order
$\ell$ on $C$ 
if and only if for any point $x\in C$ the partial derivatives 
of $g_1$ and
$g_2$ of order at most $\ell$ are equal in any local coordinate system
containing $x$.~\done 
\end{definition}

\begin{remark}\normalshape
If $C$ is the closure of an open set and $g_1=g_2$ on
$C$, then $g_1$ and $g_2$ agree to order $\ell$ on $C$ for 
all $\ell$.~\done
\end{remark}

\begin{thm}\label{flatspacerigid}\sl Let $(\am,\ame_0)$ 
be a geodesically complete
semi-Riemannian manifold of dimension at least three with constant 
curvature~0.
Let $U\subset \am$ be a connected open set with compact closure and assume
that $\pi_1(U)$ is finite. 
Further, let $\ame$ be a semi-Riemannian metric
on $\am$ and assume that $\ame=\ame_0$ to order 2 on $\am\setminus U$ and 
that the curvature tensor $\ar$ of $\am$ satisfies a one sided curvature bound
\begin{equation}
\label{curvbound2}
\ar\ge0\quad\mbox{or}\quad\ar\le0.
\end{equation}
Then $(\am,\ame)$ is isometric to $(\am,\ame_0)$  and thus 
$(\am,\ame)$ has constant sectional curvature~0.
\end{thm}

\begin{thm}\label{hypspacerigid}\sl Let $(\am,\ame_0)$ 
be a geodesically complete
Riemannian or Lorentzian manifold of dimension at least three 
which is locally isometric to $(\r^n_k(-1),g_{\r^n_k(-1)})$.  (Thus
$k=0$ (Riemannian) or $k=1$ (Lorentzian).)  Let $U\subset\am$  be a
connected open set with compact closure  and assume $\pi_1(U)$ is
finite.  Further let $\ame$ be a semi-Riemannian metric on $\am$ that
agrees with $\ame_0$ to order~2 on $\am\setminus U$ and so that the
curvature of $\ame$ satisfies a one sided curvature bound
\begin{equation}
\ar\ge-1,\quad\mbox{or}\quad\ar\le-1.		\label{curvbound1}
\end{equation}
Then $(\am,\ame)$ is isometric to $(\am,\ame_0)$.  
\end{thm}

\begin{remark}\normalshape
If the manifold $\am$ is compact, then it is possible that 
the open set $U$ of the theorems is dense in $\am$.  We look at an
example of this in the Riemannian case.  Let $(\am,\ame_0)$ be
a compact Riemannian manifold of constant sectional curvature~$K_0\le0$.  
Then the fundamental group $\pi_1(\am)$ is infinite.  Let $x\in\am$ and
let $C(x)$ be the cut locus of $x$ in $\am$.  Then $U:=\am\setminus C(x)$
is simply connected (as it is homeomorphic to a ball in $\r^n$).  Thus the
theorems say that any Riemannian metric $\ame$ on $\am$ that agrees with
$\ame_0$ to order~2 on $C(x)$ and satisfies a one sided curvature bound 
of the form (\ref{curvbound2}) or (\ref{curvbound1}),
then $\ame$ also has constant sectional curvature.~\done  
\end{remark}

\begin{pf*}{Proof of the Theorems} 
We prove Theorem~\ref{hypspacerigid},
the proof of Theorem~\ref{flatspacerigid} being similar.
The universal covering space of
$\am$ is $\r^n_k(-1)$.  Let $\pi:\r^n_k(-1)\to\am$ be the covering map.  Then 
$\pi:(\r^n_k(-1),g_{\r^n_k(-1)})\to(\am,\ame_0)$ is a local isometry.  
Let $\widehat{U}_0$ be a connected component of $\pi^{-1}[U]$.  As $\pi_1(U)$
is finite, the map $\pi|_{\widehat{U}_0}:\widehat{U}_0\to U$ is a finite 
sheeted cover.  But then using that the closure of $U$ in $\am$ is compact,
it is not hard to see that the closure of $\widehat{U}_0$ in $\r^n_k(-1)$
is compact.  Define a new semi-Riemannian metric $g$ on $\r^n_k(-1)$ by
\begin{equation}
g=\left\{\begin{array}{rl} \pi^*\ame&\mbox{\rm on $\widehat{U}_0$}\\
\pi^*\ame_0=g_{\r^n_k(-1)}&\mbox{\rm on $\r^n_k(-1)\setminus \widehat{U}_0$.}
	\end{array}\right.  \label{c2metric}
\end{equation}
As $\pi:(\r^n_k(-1),g_{\r^n_k(-1)})\to(\am,\ame_0)$ is a local
isometry and $\ame$ and $\ame_{0}$ agree to order~2 on $\partial U$ the
semi-Riemannian metric $g$ pieces together to form a $C^2$ metric on
$\r^n_k(-1)$. This semi-Riemannian metric agrees with the standard
metric $g_{\r^n_k(-1)}$ outside of the compact set 
$\mathop{\mathrm{closure}}(\widehat{U}_0)$.

It is easily checked that the proof of the rigidity result
Theorem~\ref{thm:rigidity-curved} goes through under the assumption the
metrics are $C^2$. Therefore, we are able to conclude that
$(\r^n_k(-1),g)$ is isometric to $(\r^n_k(-1),g_{\r^n_k(-1)})$.
Theorem~\ref{hypspacerigid} now follows.

The proof of Theorem~\ref{flatspacerigid} works along the same lines
except that Corollary~\ref{domain} is used to finish the proof.
\end{pf*}

\ifx\undefined\bysame
\newcommand{\bysame}{\leavevmode\hbox to3em{\hrulefill}\,}
\fi

\end{document}